\documentstyle [12pt]{article}
\let\oldphi=\phi

\let\tilde=\widetilde

\def\one#1{#1^{\raise5pt\hbox{$\scriptstyle\!\!\!\!1$}}\,{}}
\def\two#1{#1^{\raise5pt\hbox{$\scriptstyle\!\!\!\!2$}}\,{}}
\def\three#1{#1^{\raise5pt\hbox{$\scriptstyle\!\!\!\!3$}}\,{}}
\def\phi{\varphi}
\def\eps{\varepsilon}
\def\a{\alpha}
\def\b{\beta}

\def\d{\delta}
\def\D{\Delta}
\def\l{\lambda}

\def\comment#1{}
\def\abs#1{\left|#1\right|}
\def\id{\hbox{{1}\kern-.25em\hbox{\rm l}}}
\def\beq{\begin{equation}}
\def\eeq{\end{equation}}
\def\be{\begin{displaymath}}
\def\ee{\end{displaymath}}
\def\bmat{\left(\begin{array}}
\def\emat{\end{array}\right)}
\def\bds{\begin{description}}
\def\eds{\end{description}}
\def\?{(?)\marginpar{|?}}

\def\half{\frac{1}{2}}

\newtheorem{theo}{Theorem}
\newtheorem{prop}{Proposition}
\newtheorem{guess}{Conjecture}
\newtheorem{lemma}{Lemma}
\def\dd{\partial}
\def\endproof{\hfill\rule{2mm}{2mm}}
\def\D{{\cal D}}

\def\K{{\cal K}}
\def\M{{\cal M}}
\def\G{\Gamma}
\def\dfrac#1#2{\frac{\displaystyle #1}{\displaystyle #2}}
\def\abs#1{\left|#1\right|}
\def\w{\check{v}}
\def\dilog{\mathop{\rm Li}\nolimits_2}
\def\aa{\check\alpha}
\def\Lq#1#2#3{{\cal L}_q\!\left(#1;#2,#3\right)}
\def\tworow#1#2{\begin{array}{c}#1 \\ #2 \end{array} }
\newcommand{\tfrac}[2]{{\textstyle\frac{#1}{#2}}}
\renewcommand{\theequation}{\thesection.\arabic{equation}}
\newcounter{subequation}[equation]
\makeatletter
 
\expandafter\let\expandafter
\reset@font\csname reset@font\endcsname
 
\def\subeqnarray{\arraycolsep1pt
    \def\@eqnnum\stepcounter##1{\stepcounter{subequation}%
        {\reset@font\rm(\theequation\alph{subequation})}}
\jot5mm     \eqnarray}

\makeatother
\newcommand{\newsection}[1]{
\vspace{10mm}
\pagebreak[3]
\addtocounter{section}{1}
\setcounter{equation}{0}
\setcounter{subsection}{0}
\setcounter{footnote}{0}
 
\begin{flushleft}
{\Large\bf \thesection. #1}
\end{flushleft}
\nopagebreak
\medskip
\nopagebreak}
\newcommand{\newappendix}[1]{
\vspace{10mm}
\pagebreak[3]
\addtocounter{section}{1}
\setcounter{equation}{0}
 
\begin{flushleft}
{\Large\bf Appendix \thesection. #1}
\end{flushleft}
\nopagebreak
\medskip
\nopagebreak}
\newfont{\bbd}{msbm10 scaled\magstep1} 
\def\C{\hbox{\bbd C}}                  
\def\R{\hbox{\bbd R}}                  
\def\Z{\hbox{\bbd Z}}
\begin{document}
\begin{flushright}
 \sf q-alg/9602023 \\
 \sf 25/3/1996
\end{flushright}
\vskip 1cm
\begin{center} \LARGE\bf
 Separation of variables \\
 for $A_2$ Ruijsenaars model \\
and new integral representation \\
 for $A_2$ Macdonald polynomials
\end{center}
\vskip 0.5cm
\begin{center}
V.B.~Kuznetsov
\footnote{Part of this work was done while the authors 
were visiting the Institute of Mathematical Modelling,
Technical University of Denmark, Lyngby, Denmark.}${}^,$
\footnote{On leave from Department of Mathematical and Computational
Physics, Institute of Physics, St.~Petersburg University, St.~Petersburg
198904, Russia.} \\
\vskip 0.2cm
Department of Applied Mathematical Studies, University of Leeds, \\
Leeds LS2 9JT, UK\\
\vskip 0.3cm
E.K.~Sklyanin ${}^1$ \\
\vskip 0.2cm
Steklov Mathematical Institute, Fontanka 27, St.Petersburg 191011, Russia
\end{center}
\vskip 0.8cm
\begin{center}
\bf Abstract
\end{center}
Using the Baker-Akhiezer function technique we construct
a separation of variables for the classical trigonometric 3-particle 
Ruijsenaars model (relativistic generalization of Calogero-Moser-Sutherland 
model). In the quantum case, an integral operator $M$ is constructed 
from the Askey-Wilson contour integral. The operator $M$ transforms
the eigenfunctions of the commuting Hamiltonians (Macdonald polynomials for 
the root sytem $A_2$) into the factorized form $S(y_1)S(y_2)$ where
$S(y)$ is a Laurent polynomial of one variable expressed in terms
of the ${}_3\oldphi_2(y)$ basic hypergeometric series. The inversion of $M$
produces a new integral representation for the $A_2$ Macdonald polynomials.
We also present some results and conjectures for general $n$-particle case.
\vskip 0.5cm
\begin{flushleft}
{\it Submitted to \rm J.~Phys.~A: Math.~Gen.}
\end{flushleft}

\newpage
\newsection{Introduction}

The Separation of Variables (SoV) is an approach to  quantum
integrable systems which can be briefly formulated as follows (for a more
detailed discussion see the survey \cite{Skl-rev}).

Given a quantum-mechanical system of $n$ degrees of freedom possessing
$n$ commuting Hamiltonians
\beq
 [H_j,H_k]=0, \qquad j,k=1,2,\ldots,n
\label{eq:comm-H}
\eeq
one tries to find an operator $M$ sending any common eigenvector 
$P_{\l}$ of the Hamiltonians 
\beq
 H_jP_{\l}=h_j P_{\l}
\label{eq:spec-H}
\eeq
labelled by the quantum numbers
$\l=\{\l_1,\ldots,\l_n\}$ into the product
\beq
 M: P_{\l}\rightarrow \prod_{j=1}^n S_{\l;j}(y_j)
\label{eq:def-SoV}
\eeq
of functions $S_{\l;j}(y_j)$ of one variable each. The original
multi-dimensional eigenvalue problem (\ref{eq:spec-H})
is transformed respectively into a set of simpler one-dimensional spectral 
problems (separated equations)
\beq
 {\cal D}_j\left(y_j,\frac{\partial}{\partial y_j};
h_1,\ldots,h_n\right)S_{\l;j}(y_j)=0
\label{eq:sep-eq-gen}
\eeq
where ${\cal D}_j$ are usually some differential
or finite-difference operators in variable $y_j$ depending on the
spectral parameters $h_k$. In the context of the classical Hamiltonian
mechanics the above construction corresponds precisely to the 
standard definition of SoV in the Hamilton-Jacobi equation.

The advent of the Inverse Scattering Method gave new life to SoV 
providing it with the interpretation of the separated coordinates $y_j$
(in the classical case) as the poles of the Baker-Akhiezer function
(properly normalized eigenvector of the corresponding Lax matrix).
The unsolved question is, however, how to choose a correct normalization
of B-A function to obtain SoV for a given Lax matrix. Nevertheless, as
an heuristic recipe, the above idea has proved to be quite efficient
and allowed to find SoV for a few new classes of classical integrables
systems. In particular, SoV was found for the systems arising from the
$r$-matrices satisfying the classical Yang-Baxter equation in case of
$A_{n-1}$ ($sl_n$) Lie algebra. In the cases $n=2$ and $n=3$ 
the construction of SoV has been successfully transferred to the quantum case
(see \cite{Skl-rev} and references therein).

Pursuing the goal to extend the applicability of the B-A function recipe,
in our previous paper \cite{Kuz-Skl} we have studied the $A_2$
Calogero-Sutherland model which does not fall into the previously studied
cases since it posesses a dynamical (non-numeric)
$r$-matrix \cite{dynam-r}. In the quantum case, our construction of SoV
has produced a new integral representation for the eigenfunctions of the
$A_2$ C-S Hamiltonians (known as Jack polynomials) in terms of
${}_3F_2$ hypergeometric polynomials.

In the present paper we generalize the results of \cite{Kuz-Skl} to the 
3-particle Ruijsenaars model \cite{Ruijs}
which is a relativistic analog of the C-S model.
The corresponding eigenfunctions (Macdonald polynomials 
\cite{Macd-book,Macd-poly}) are $q$-analogs of Jack polynomials. No surprize
that the corresponding separated functions 
are Laurent polynomials expressed in terms of
${}_3\oldphi_2$ basic hypergeometric series. 
We present also some results and conjectures for the general
$n$-particle problem, for instance, we connect the $A_{n-1}$ type
basic hypergeometric separation polynomials $S_{\l}(y)$ to 
a terminated case of the $\oldphi_D$ type $q$-Lauricella function
of $n-1$ variables.

The paper is organized as follows. In Section 2 we describe the classical
Ruijsenaars model and, using B-A function technique, construct a SoV.
Though the results of this section are not used directly in what follows,
they provide a useful background for subsequent treatment of the quantum
case. 
In Section 3 the standard facts concerning the quantum Ruijsenaars model
and Macdonald polynomials are collected.
In Section 4, after introducing the quantum Hamiltonians and
Macdonald polynomials, we describe the integral operator $M$ performing
a SoV and formulate the main theorem whose proof takes the rest of the
section and part of the next one. 
The main part of the proof is contained in Section 4 where the
properties of the operator $M$ are studied, whereas in Section 5 the 
results concerning the separated equation 
(certain 3-rd order $q$-difference equation and its $n$-th order
generalization), as well as  its polynomial solutions, are collected.
The main technical tool allowing us to study the operator $M$ is the famous
Askey-Wilson integral identity (\ref{eq:AW-int}).

Generally, SoV is aimed to simplify the multidimensional spectral problem
by reducing it to a series of one-dimensional ones. In case of the 
Calogero-Sutherland and Ruijsenaars models, however, the spectrum and
eigenfunctions are well known  and studied by independent means.
The main benefit of SoV in application to these models is rather
producing new relations between special functions. In particular,
inverting the operator $M$ one obtains a new integral representation
for $A_2$ Macdonald polynomials in terms of ${}_3\oldphi_2$ basic 
hypergeometric functions, which is done in the end of Section 4.
In section 6 we discuss the obtained results and the possibility of their
generalization to $A_n$, $n>2$ case.
Two Appendices, A and B, contain, respectively, a collection of necessary
formulas from $q$-analysis and some auxiliary results concerning 
operator $M$.

During our work we enjoyed the hospitality of the Research Institute
for Mathematical Sciences (Kyoto, Japan) and the Institute of Mathematical
Modelling (Lyngby, Denmark). We express our gratitude to our hosts
Prof.\ T.~Miwa and Prof.\ P.~L.~Christiansen.
VBK acknowledges supports from the grant of Forskerakademiet, Aarhus,
while in the Technical University of Denmark and grant from CRM-ISM
as a post-doctoral fellow in the University of Montreal.

\newsection{Classical Ruijsenaars model}

In the spirit of $q$-analysis,
we prefer to use exponentiated canonical coordinates and momenta.

{\bf Definition 1.} {\it The variables $(X_j,x_j)$ $j=1,\ldots,n$ 
on a $2n$-dimensional symplectic manifold form
a {\rm Weyl canonical system} if they possess the Poisson brackets
\beq
 \{X_j,X_k\}=\{x_j,x_k\}=0, \quad \{X_j,x_k\}=-iX_jx_k\d_{jk}, 
 \qquad j,k=1,\ldots,n
\label{eq:def-Weyl-can}
\eeq
or, equivalently, the symplectic form $\omega$ is expressed as
$\omega=i\sum_j d\ln X_j\wedge d\ln x_j=
d(i\sum_j\ln X_j d\ln x_j)$.
}

The $n$-particle ($A_{n-1}$) trigonometric Ruijsenaars model \cite{Ruijs}
is formulated in terms of the Weyl canonical  system $(T_j,t_j)$ where
$\abs{t_j}=1$, $T_j\in\R$ $(j=1,2,\ldots,n)$.
The Hamiltonians $H_i$ are defined as 
\beq
 H_i=\sum_{J\subset\{1,\ldots,n\} \atop \left|J\right|=i}
  \left(\prod_{j\in J \atop k\in \{1,\ldots,n\}\setminus J} v_{jk} \right)
  \left(\prod_{j\in J}T_j \right)\,, \qquad i=1,\ldots,n,
\label{eq:def-H}
\eeq
where
\beq
 v_{jk}=\frac{\ell^{-\half}t_j-\ell^{\half}t_k}{t_j-t_k}\,, \qquad
 \ell\in(1,\infty)\,.
\label{eq:def-w}
\eeq

\begin{prop} {\bf\cite{Ruijs,vD}}

The Hamiltonians $H_j$ Poisson commute.
\beq
 \{H_j,H_k\}=0, \qquad j,k=1,\ldots,n.
\label{eq:comm-H-cl}
\eeq
\end{prop}

Define the Lax matrix ($L$ operator) by the formula
\beq
   L(u)=D(u)E(u)
\label{eq:L=DE}
\eeq
where
\beq
 D_{jk}=\frac{(\ell-1)(1-\ell^nu)}{2\ell^{\frac{n+1}{2}}(1-u)}
\left(\prod_{i\neq j}v_{ji}\right)T_j\d_{jk},
\label{eq:def-D-gen}
\eeq
\beq
  E_{jk}=\frac{1+\ell^nu}{1-\ell^nu}-\frac{t_j+\ell t_k}{t_j-\ell t_k}.
\label{eq:def-E-gen}
\eeq

\begin{prop} {\bf\cite{Ruijs}}

The characteristic polynomial of the matrix L(u) (\ref{eq:L=DE})
generates the Hamiltonians (\ref{eq:def-H})
\begin{eqnarray}
\lefteqn{
 (-1)^n\ell^{\frac{n(n-1)}{2}}(1-\ell^nu)(1-u)^n\det(z-L(u))
} \nonumber \\
&&=\sum_{k=0}^n(-1)^k\ell^{\frac{n-1}{2}k}(1-\ell^ku)(1-u)^k(1-\ell^nu)^{n-k}
H_{n-k}z^k
\label{eq:char-poly-N}
\end{eqnarray}
where we assume $H_0\equiv1$.
\end{prop}

In the 3-particle ($A_2$) case which we consider henceforth we have,
respectively,
\begin{subeqnarray}
 H_1&=&v_{12}v_{13}T_1+v_{21}v_{23}T_2+v_{31}v_{32}T_3, \\
 H_2&=&v_{13}v_{23}T_1T_2+v_{12}v_{32}T_1T_3+v_{21}v_{31}T_2T_3, \\
 H_3&=&T_1T_2T_3,
\label{eq:def-H-3}
\end{subeqnarray}

\beq
 D=\frac{(\ell-1)(1-\ell^3u)}{2\ell^2(1-u)}\,
  \mbox{\rm diag}
\left\{v_{12}v_{13}T_1,v_{21}v_{23}T_2,v_{31}v_{32}T_3\right\},
\label{eq:def-D-3}
\eeq
\beq
  E_{jk}=\frac{1+\ell^3u}{1-\ell^3u}-\frac{t_j+\ell t_k}{t_j-\ell t_k},
\label{eq:def-E-3}
\eeq
and
\begin{eqnarray}
 \lefteqn{\kern-2cm\ell^3(1-u)^2\det(z-L(u))=
z^3\ell^3(1-u)^2-z^2\ell^2(1-u)(1-\ell^2u)H_1 } \nonumber \\
&&   +z\ell(1-\ell u)(1-\ell^3 u)H_2-(1-\ell^3u)^2H_3.
\label{eq:char-pol-cl-3}
\end{eqnarray}

To find a SoV for the Ruijsenaars system we use the recipe discussed in
the Introduction and choose for the separated coordinates $y_j$ the poles
upon $u$  of the Baker-Akhiezer function $\psi(u)$ (an eigenvector
of $L(u)$) normalized 
by the condition that its 3-rd component $\psi_3(u)$ is constant. The
canonically conjugated (in the Weyl sense) variables $Y_j$ are chosen as
the eigenvalues of $L(y_j)$. For the detailed discussion of the B-A function
recipe see \cite{Skl-rev} though the construction described below is quite
self-contained.

Define two functions $A_1(u)$ and $A_2(u)$ by the formulas
\beq
 A_k(u):=L_{kk}-\frac{L_{3k}L_{k,3-k}}{L_{3,3-k}}=T_k\a_k(u), \qquad k=1,2
\label{eq:def-A}
\eeq
\beq
 \a_k(u):=\frac{(1-\ell^3u)(\ell t_3u-t_{3-k})(t_k-\ell t_3)}%
{\ell(1-u)(\ell^2 t_3u-t_{3-k})(\ell t_k-t_3)}, \qquad k=1,2\,.
\label{eq:def-alpha-cl}
\eeq

The separated variables $y_j$ are defined from the equation
\beq
   A_1(y)=A_2(y).
\label{eq:def-y}
\eeq

It is easy to see that (\ref{eq:def-y}) has 3 solutions one of which 
$y=\ell^{-3}$ we ignore since it is a constant. The remaining two roots
we denote $y_1$ and $y_2$. From the easily verified invariance of
$\a_1(u)/\a_2(u)$ under the transformation 
$u\mapsto u^{-1}t_1t_2t_3^{-2}\ell^{-\frac32}$
it follows that
\beq
 y_1y_2=\frac{t_1t_2}{t_3^2\ell^3}.
\label{eq:cl-constraint}
\eeq

The conjugated variables $Y_j$ are defined as
\beq
 Y_j=A_1(y_j)=A_2(y_j), \qquad j=1,2.
\label{eq:def-Y}
\eeq

Equivalently, the four variables $Y_1,Y_2,y_1,y_2$ are defined through
four equations
\beq
 Y_j=T_k\a_k(y_j), \qquad j,k\in\{1,2\}.
\label{eq:four-eq-Yy}
\eeq
\begin{theo}
The variables $Y_j$, $y_j$ satisfy the separated equations
\begin{eqnarray}
\lefteqn{
Y_j^3\ell^3(1-y_j)^2-Y_j^2\ell^2(1-y_j)(1-\ell^2y_j)H_1 }\nonumber \\
 &&  +Y_j\ell(1-\ell y_j)(1-\ell^3 y_j)H_2-(1-\ell^3y_j)^2H_3=0,
\qquad j=1,2
\label{eq:char-eq-cl}
\end{eqnarray}
which, by virtue of (\ref{eq:char-pol-cl-3}), imply that 
$\det(Y_j-L(y_j))=0$.
\label{th:cl-char-eq}
\end{theo}
{\bf Proof.}
Substitute into (\ref{eq:char-eq-cl}) the expressions (\ref{eq:def-H-3})
for $H_j$ and split the left-hand-side of (\ref{eq:char-eq-cl}) into
two terms
\beq
 T_3Z_1+Y_jZ_2=0
\label{eq:Z1+Z2}
\eeq
where
\begin{subeqnarray}
\lefteqn{ \kern-5mm Z_1=-(1-y_j)(1-\ell^2y_j)\ell^2v_{31}v_{32}Y_j^2
}\nonumber \\
&&\kern-5mm +(1-\ell y_j)(1-\ell^3y_j)\ell
(v_{12}v_{32}T_1Y_j+v_{21}v_{31}T_2Y_j) 
-(1-\ell^3y_j)^2T_1T_2, \\
Z_2&=&(1-y_j)^2\ell^3Y_j^2 
-(1-y_j)(1-\ell^2y_j)\ell^2(v_{12}v_{13}T_1Y_j+v_{21}v_{23}T_2Y_j)
\nonumber \\
&&+(1-\ell y_j)(1-\ell^3y_j)v_{13}v_{23}T_1T_2.
\label{eq:def-Z12}
\end{subeqnarray}

To prove (\ref{eq:char-eq-cl}) it is sufficient to show that $Z_1=Z_2=0$.
Replacing $Y_j$ in (\ref{eq:def-Z12}) by $T_1\a_1(y_j)$ or $T_2\a_2(y_j)$
in such a way that the factor $T_1T_2$ could be cancelled from
$Z_{1,2}$ we obtain that $Z_{1,2}=0$ follows from two
algebraic identities for $\a_{1,2}$
\begin{subeqnarray}
\lefteqn{
-(1-y)(1-\ell^2y)\ell^2v_{31}v_{32}\a_1(y)\a_2(y)
} \nonumber \\
&&+(1-\ell y)(1-\ell^3y)\ell
\bigl(v_{12}v_{32}\a_2(y)+v_{21}v_{31}\a_1(y)\bigr)
-(1-\ell^3y)^2=0, \\
\lefteqn{
(1-y)^2\ell^3\a_1(y)\a_2(y)
-(1-y)(1-\ell^2y)\ell^2
\bigl(v_{12}v_{13}\a_2(y)+v_{21}v_{23}\a_1(y) \bigr)
} \nonumber \\
&& +(1-\ell y)(1-\ell^3y)\ell v_{13}v_{23}
=0,
\label{eq:alpha-cl-id}
\end{subeqnarray}
which are verified directly.
\endproof

The third pair of separated variables is defined as
\beq
  x:=t_3, \qquad X:=T_1T_2T_3,
\label{eq:def-Xx}
\eeq
the corresponding separated equation being
\beq
 X-H_3=0.
\label{sep-eq-X}
\eeq

\begin{theo}
The variables $(X,Y_1,Y_2;x,y_1,y_2)$ form a Weyl canonical system in the
sense of the definition 1.
\end{theo}
{\bf Proof.}
Let us introduce new variables:
\begin{subeqnarray}
 &t_+=t_1^{1/2}t_2^{1/2}t_3^{-1}, \qquad & t_-=t_1^{1/2}t_2^{-1/2},   \\
 & T_+=T_1T_2, \qquad & T_-=T_1T_2^{-1},
\label{eq:def-tpm}
\end{subeqnarray}
and also
\begin{subeqnarray}
 &y_+=y_1^{1/2}y_2^{1/2}, \qquad & y_-=y_1^{1/2}y_2^{-1/2},   \\
 &Y_+=Y_1Y_2, \qquad & Y_-=Y_1Y_2^{-1}.
\label{eq:def-ypm}
\end{subeqnarray}

Obviously, $(X,T_-,T_+;x,t_-,t_+)$ is also a Weyl canonical system.
Note that 
\beq
 y_+=t_+\ell^{-\frac{3}{2}}
\label{eq:yp=tp}
\eeq
because of (\ref{eq:cl-constraint}). Note also that from 
(\ref{eq:def-alpha-cl}) it follows that $Y_\pm$, $y_\pm$ depend only on 
$T_\pm$, $t_\pm$ and do not contain $X$, $x$.

It remains to show that the transformation from 
$(T_-,T_+;t_-,t_+)$ to $(Y_-,Y_+;$ $y_-,y_+)$ is canonical that is
$(Y_-,Y_+;y_-,y_+)$ is again a Weyl canonical system. To this end,
it suffices to construct the generating function $F(Y_+,y_-;t_+,t_-)$
of the canonical transformation such that \cite{Arnold}
\beq
{i}\ln T_\pm=t_\pm\frac{\dd F}{\dd t_\pm},\qquad
{i}\ln Y_-=-y_-\frac{\dd F}{\dd y_-},\qquad
{i}\ln y_+=Y_+\frac{\dd F}{\dd Y_+}. 
\label{eq:gen-func}
\eeq
and $d(F-{i}\ln Y_+\ln y_+)=
{i}(\ln T_-d\ln t_-+\ln T_+d\ln t_+)-
{i}(\ln Y_-d\ln y_-+\ln Y_+d\ln y_+)$. 
Recalling the definition of the Euler dilogarithm \cite{dilog}
\beq
 \dilog(z):=-\int_0^z\frac{dt}{t}\ln(1-t)=\sum_{k=1}^\infty\frac{z^k}{k^2}
\label{eq:dilog}
\eeq
and introducing the notation 
\beq
 {\cal L}(\nu;x,y):=\dilog(\nu xy)+\dilog(\nu xy^{-1})
+\dilog(\nu x^{-1}y)+\dilog(\nu x^{-1}y^{-1})
\label{eq:def-LL}
\eeq
we define $F:=i\ln Y_+\ln(\ell^{-\frac32}t_+)+\tilde F$,
\begin{eqnarray}
 \tilde F&:=&{i}\bigl({\cal L}(\ell^{-\half};y_-,t_-)
+{\cal L}(\ell^{-1};t_+,t_-)-{\cal L}(\ell^{-\frac32};t_+,y_-) \nonumber \\
&&-\dilog(t_-^2)-\dilog(t^{-2}_-)\bigr).
\label{eq:def-F}
\end{eqnarray}

It is a matter of direct calculation to verify, using (\ref{eq:four-eq-Yy})
and (\ref{eq:yp=tp}), that $F$ satisfies (\ref{eq:gen-func}).
\endproof

The identities (\ref{eq:char-eq-cl}) and 
(\ref{sep-eq-X}) and canonicity of the variables $(X,Y_1,Y_2;x,y_1,$ $y_2)$
established above provide, by definition \cite{Skl-rev}, a SoV for the
$A_2$ Ruijsenaars system.

\newsection{Quantization}
We collect here the standard facts concerning the quantum
$n$-particle ($A_{n-1}$) Ruijsenaars model \cite{Ruijs,vD} and the 
corresponding Macdonald polynomials \cite{Macd-book,Macd-poly}.

Throughout the paper $\Z$ stands for the set of integers, the notations
$\Z_{\geq0}$ and $\Z_{\leq0}$ are self-evident.

The quantum Ruijsenaars model is described in terms of the multiplication
and shift
operators, resp.\ $t_j$ and $T_j$ ($j=1,\ldots,n$) acting on functions
of $t_j$
\beq
 (t_jf)(\vec t):=t_jf(\vec t), \qquad
 (T_jf)(\vec t):=f(\ldots,qt_j,\ldots)
\label{eq:def-T}
\eeq
(we do not make distinction between variables and operators $t_j$).
Here $q$ is the quantum deformation parameter related to the Planck constant
$\hbar>0$ as
\beq
 q=e^{-\hbar}, \qquad q\in(0,1).
\label{def-q}
\eeq

The operators $T_j$, $t_j$ satisfy the Weyl commutation relations
\beq
 [T_j,T_k]=[t_j,t_k]=0,\qquad
 T_jt_k=\left\{\begin{array}{ll}
         qt_kT_j, & j=k \\
          t_kT_j, & j\neq k \end{array}\right.
\label{eq:comm-Pt}
\eeq
which produce the Poisson brackets (\ref{eq:def-Weyl-can})
in the classical limit $\hbar\rightarrow0$ by the standard correspondence
rule $[,]=-i\hbar\{,\}+O(\hbar^2)$.

The commuting quantum Hamiltonians $H_j$
\beq
 [H_j,H_k]=0, \qquad j,k=1,\ldots,n
\label{eq:comm-S-q}
\eeq
are given by the same formulas (\ref{eq:def-H}) as in the classical case
with the fixed operator ordering ($T_j$ to the right). We assume that
\beq
 \ell=q^{-g}=e^{g\hbar}, \qquad g>0, \quad \ell\in(1,\infty)
\label{eq:def-ell-g}
\eeq
(note that both in the classical and nonrelativistic limits 
$\hbar\rightarrow0$, $q=e^{-\hbar}\rightarrow 1$ but in the
classical limit $g\rightarrow\infty$, $\ell={\rm const}$ whereas
in the nonrelativistic limit $g={\rm const}$, $\ell\rightarrow1$).

The operators $H_k$ leave invariant the space ${\rm Sym}(t_1,\ldots,t_n)$ of
symmetric Laurent polynomials in variables $t_j$. 
A basis in ${\rm Sym}(t_1,\ldots,t_n)$  is given by the
monomial symmetric functions $m_{\l}$
labelled by the sequences 
$\l=\{\l_1\leq \l_2\leq\ldots\leq \l_n\}$ of integers $\l_j\in\Z$
(dominant weights)
and expressed as 
$m_{\l}=\sum t_1^{\nu_1}\ldots t_n^{\nu_n}$ where the sum is taken over
all distinct permutations $\nu$ of $\l$. 

Denote $\abs{\l}\equiv\sum_{j=1}^n \l_j$.
The dominant order on the dominant weights $\l$ is defined as
\beq
 \l'\preceq \l \quad \Longleftrightarrow \quad
\left\{ 
\abs{\l'}=\abs{\l};
\quad
\sum_{j=k}^n \l_j^\prime\leq\sum_{j=k}^n \l_j, \quad
k=2,\ldots,n\right\}.
\label{eq:dom-ord}
\eeq

The Macdonald polynomials 
$P_{\l}^{(\ell;q)}\in{\rm Sym}(t_1,\ldots,t_n)$ are
uniquely defined as 
joint eigenvectors of $H_k$ in ${\rm Sym}(t_1,\ldots,t_n)$
\beq
 H_k P_{\l}^{(\ell;q)}=h_kP_{\l}^{(\ell;q)}
\label{eq:eigv-S}
\eeq
labelled by the dominant
weight $\l$ and normalized by the condition
\beq
P_{\l}^{(\ell;q)}=
\sum_{\l'\preceq \l}\kappa_{\l\l'}m_{\l'}, \qquad
 \kappa_{\l\l}=1.
\label{eq:def-Macd-poly}
\eeq

The corresponding eigenvalues $h_k$ are
\beq
 h_k=\sum_{j_1<\ldots<j_k}\mu_{j_1}\ldots \mu_{j_k}, \qquad
 \mu_j=q^{\l_j}\ell^{\frac{n+1}{2}-j}.
\label{eq:def-hk-gen}
\eeq

Note that our parameter $\ell$ and parameter $t$ used in 
\cite{Macd-book,Macd-poly} relate as $\ell=t^{-1}$.

The polynomials $P_{\l}^{(\ell;q)}$ are orthogonal
\beq
 \frac{1}{(2\pi i)^n}\oint\limits_{\abs{t_1}=1}\frac{dt_1}{t_1}\ldots
\oint\limits_{\abs{t_n}=1}\frac{dt_n}{t_n}
\bar P_{\l}^{(\ell;q)}(\vec t)P_{\l'}^{(\ell;q)}(\vec t)
\Delta(\vec t)=0, \qquad
\l\neq \l'
\label{eq:orthog-Macd}
\eeq
with respect to the weight
\beq
 \Delta(t_1,\ldots,t_n)=\prod_{j\neq k}
\frac{(t_jt_k^{-1};q)_\infty}{(\ell^{-1}t_jt_k^{-1};q)_\infty}
\label{eq:def-Delta}
\eeq
(see (\ref{eq:Poch-inf}) for the notation).

In the limit $\hbar\rightarrow0$, $g={\rm const}$ 
the appropriate linear combinations of $H_k$
produce the Hamiltonians of the nonrelativistic Calogero-Sutherland model,
and the Macdonald polynomials go over into the Jack polynomials,
see \cite{Kuz-Skl}.

In the present paper we consider only the simplest nontrivial case $n=3$.

The Hamiltonians $H_k$ being given by (\ref{eq:def-H-3}),
the formulas (\ref{eq:def-hk-gen}) produce, respectively, 
\beq
 h_1=\ell q^{\l_1}+q^{\l_2}+\ell^{-1}q^{\l_3}, \quad
 h_2=\ell q^{\l_1+\l_2}+q^{\l_1+\l_3}+\ell^{-1}q^{\l_2+\l_3}, \quad
 h_3=q^{\abs{\l}}
\label{eq:def-hk-3}
\eeq
for their eigenvalues labelled by the ordered triplets
$ \{\l_1\leq \l_2\leq \l_3 \}\in {\Z}^3$.

For instance,
\be
 m_{000}=1, \quad
 m_{001}=t_1+t_2+t_3, \quad
 m_{011}=t_1t_2+t_1t_3+t_2t_3, \quad
 m_{002}=t_1^2+t_2^2+t_3^2,
\ee
\be
 m_{111}=t_1t_2t_3, \quad
 m_{012}=t_1t_2^2+t_1^2t_2+t_1t_3^2+t_1^2t_3+t_2t_3^2+t_2^2t_3,
\ee
\be
 m_{112}=t_1^2t_2t_3+t_1t_2^2t_3+t_1t_2t_3^2, \quad
 m_{022}=t_1^2t_2^2+t_1^2t_3^2+t_2^2t_3^2, \quad
 m_{003}=t_1^3+t_2^3+t_3^3.
\ee

\be
 P^{(\ell;q)}_{000}=m_{000}, \quad
 P^{(\ell;q)}_{001}=m_{001}, \quad
 P^{(\ell;q)}_{011}=m_{011}, \quad
 P^{(\ell;q)}_{002}=m_{002}+\tfrac{(1-\ell)(1+q)}{q-\ell}m_{011},
\ee
\be
 P^{(\ell;q)}_{111}=m_{111}, \quad
 P^{(\ell;q)}_{012}=m_{012}+
\tfrac{(1-\ell)(q(2+\ell)+1+2\ell)}{q-\ell^2}m_{111},
\ee
\be
 P^{(\ell;q)}_{112}=m_{112}, \quad
 P^{(\ell;q)}_{022}=m_{022}+\tfrac{(1-\ell)(1+q)}{q-\ell}m_{112},
\ee
\be
 P^{(\ell;q)}_{003}=m_{003}+\tfrac{(1-\ell)(1+q+q^2)}{q^2-\ell}m_{012}
+\tfrac{(1-\ell)^2(1+q)(1+q+q^2)}{(q-\ell)(q^2-\ell)}m_{111}.
\ee

\newsection{Operator $M$}

We are now going to describe the integral operator $M$ (\ref{eq:def-SoV})
producing the SoV. Generally speaking, the kernel $\M$ of $M$ should
depend on 6 variables: $\M(x,y_1,y_2\mid t_1,\allowbreak t_2,t_3)$.
However, by analogy with the classical case (section 2) and the
nonrelativistic limit \cite{Kuz-Skl}, it is natural to assume that $\M$
contains two $\d$-functions corresponding to the constraints $x=t_3$
(\ref{eq:def-Xx}) and, respectively, (\ref{eq:cl-constraint}).
There remains thus only one integration in $M$. Again by analogy with
the previously studied cases, the kernel $\M$ is most conveniently described
in terms of the variables $t_\pm$ (\ref{eq:def-tpm}a) and
$y_\pm$ (\ref{eq:def-ypm}a).

So, let us introduce the operator $M$
\begin{eqnarray}
\lefteqn{
 M: \Psi(t_1,t_2,t_3)\rightarrow\Phi(x,y_1,y_2) }\nonumber \\
&&\kern-12mm=\frac{1}{2\pi i}
\int\limits_{\G_{g,2g}^{t_+,y_-}} \!\!\frac{dt_-}{t_-}
\M\bigl((y_1y_2)^{\frac12},(y_1/y_2)^{\frac12}\!\bigm|\! t_-\bigr)
\Psi\bigl(\ell^{\frac32}x(y_1y_2)^{\frac12}t_-,
     \ell^{\frac32}x(y_1y_2)^{\frac12}t_-^{-1},x\bigr)
\label{eq:def-M}
\end{eqnarray}
with the kernel
\beq
\M(y_+,y_-\mid t_-)=
\frac{(1-q)(q;q)_\infty^2(t_-^2,t_-^{-2};q)_\infty
\,\Lq{\ell^{-\frac32}}{y_-}{y_+\ell^{\frac32}} }%
{2B_q(g,2g)
\,\Lq{\ell^{-\half}}{y_-}{t_-}\Lq{\ell^{-1}}{t_-}{y_+\ell^{\frac32}} }
\label{eq:def-ker-M}
\eeq
where the notation (\ref{eq:def-q-Beta}) and (\ref{eq:def-Lq}) is used.
For the definition of the cycle $\G_{g,2g}^{t_+,y_-}$ which depends on 
$g$, $y_{1,2}$ see (\ref{eq:def-Kab}) and (\ref{eq:def-contour}).

{\bf Remark.} In the classical limit, as $q\rightarrow1$, 
$\ell={\rm const}$, using (\ref{eq:asympt-dilog}) and
$\ln{\cal L}_q(\nu;\allowbreak x,\allowbreak y)
\sim-\hbar^{-1}{\cal L}(\nu;x,y)$ 
one obtains that
the asymptotics $\ln\M\sim-i\hbar^{-1}\tilde F$
of the kernel $\M$ is determined by the regular part
$\tilde F$ (\ref{eq:def-F}) of the generating function of the 
canonical transformation producing classical SoV.
As for the nonrelativistic limit, $\hbar\rightarrow0$, $g={\rm const}$,
the easiest way to reproduce the results of \cite{Kuz-Skl} is to compare
the action of the operators $M$ and its nonrelativistic analog on 
polynomials, see theorem \ref{th:act-M-basis-t}.

Now we are in a position to formulate our main result.
\begin{theo}
The operator $M$ (\ref{eq:def-M}) transforms any $A_2$ Macdonald
poly\-no\-mi\-al \break
$P_{\l}^{(\ell;q)}(t_1, t_2, t_3)$ into the product
\beq
 M: P_{\l}^{(\ell;q)}(t_1, t_2, t_3)\rightarrow
         c_{\l}x^{\abs{\l}}
S_{\l}^{(\ell;q)}(y_1)S_{\l}^{(\ell;q)}(y_2)
\label{eq:act-M-Macd-3}
\eeq
of functions of one variable only, where the Laurent polynomials
$S_{\l_1\l_2\l_3}^{(\ell;q)}$
\beq
 S_{\l}^{(\ell;q)}(y)=\sum_{k=\l_1}^{\l_3} \chi_{\l,k}^{(\ell;q)}y^k
\label{eq:pn-poly}
\eeq
are expressed in terms of the basic hypergeometric series
(\ref{eq:def-bas-hgs})
\beq
 S_{\l}^{(\ell;q)}(y)=y^{\l_1}(y;q)_{1-3g}\,
   {}_3\/\oldphi_2\left[\begin{array}{c}
     \ell^3q^{1-\l_{31}},\ell^2q^{1-\l_{21}},\ell q\\
     \ell^2q^{1-\l_{31}},\ell q^{1-\l_{21}} \end{array};q,y\right],
\label{eq:def-sep-poly-3}
\eeq
where $\l_{jk}\equiv \l_j-\l_k$. The coefficients $\chi_{\l,k}^{(\ell;q)}$
are given by
\beq
 \chi_{\l,k}^{(\ell;q)}=(q\ell^3)^{k-\l_1}
\frac{(q^{-1}\ell^{-3};q)_{k-\l_1}}{(q;q)_{k-\l_1}}
{}_4\/\oldphi_3\left[\begin{array}{c}
     q^{\l_1-k},\ell^3q^{1-\l_{31}},\ell^2q^{1-\l_{21}},\ell q\\
 \ell^3
q^{\l_1-k+2},\ell^2q^{1-\l_{31}},\ell q^{1-\l_{21}} \end{array};q,q\right].
\label{eq:def-chi-3}
\eeq

The normalization coefficient $c_{\l}$ equals
\beq
  c_{\l}=\ell^{4\l_1-\l_2}
\frac{(\ell^{-2};q)_{\l_{31}}(\ell^{-2};q)_{\l_{32}}(\ell^{-1};q)_{\l_{21}}}%
{     (\ell^{-3};q)_{\l_{31}}(\ell^{-1};q)_{\l_{32}}(\ell^{-2};q)_{\l_{21}}}.
\label{eq:def-c-3}
\eeq
\label{th:main-th}
\end{theo}

The proof of the above result will occupy the rest of this section and
a part of the next one. Our proof parallels the similar one for
the nonrelativistic Calogero-Sutherland model \cite{Kuz-Skl}.

We begin with proving the factorization (\ref{eq:act-M-Macd-3}) of
$MP_{\l}^{(\ell;q)}$.
The first step is to show that the image $MP_{\l}^{(\ell;q)}$
satisfies certain $q$-difference equations in $x$, $y_1$, $y_2$.
Let us introduce the operators $Y_j$ ($j=1,2$) acting on functions
of $y_k$ as (cf.\ (\ref{eq:def-T}))
\beq
 (Y_jf)(\vec y)=f(\ldots,qy_j,\ldots).
\label{eq:act-Yy}
\eeq

Using $y_\pm=(y_1y_2^{\pm1})^{1/2}$ (\ref{eq:def-ypm}a) one can write also
\beq
 (Y_1f)(y_+,y_-)=f(q^{\half}y_+,q^{\half}y_-), \qquad
 (Y_2f)(y_+,y_-)=f(q^{\half}y_+,q^{-\half}y_-).
\label{eq:act-Y-ypm}
\eeq

Similarly,
\beq
 (T_1f)(t_+,t_-)=f(q^{\half}t_+,q^{\half}t_-), \qquad
 (T_2f)(t_+,t_-)=f(q^{\half}t_+,q^{-\half}t_-).
\label{eq:act-T-tpm}
\eeq

We define also the operator $X$ as $X(f)(x)=f(qx)$.

Let us introduce the operator expression ${\cal D}$
\begin{eqnarray}
\lefteqn{\kern-9mm
{\cal D}(u,z;{\cal H}_1,{\cal H}_2,{\cal H}_3):=
(1-qu)(1-q^2u)\ell^3z^3
-(1-qu)(1-q^2\ell^2u)\ell^2z^2{\cal H}_1 }\nonumber \\
&&+(1-q\ell u)(1-q^2\ell^3u)\ell z {\cal H}_2
-(1-q\ell^3u)(1-q^2\ell^3u){\cal H}_3
\label{eq:def-D-op}
\end{eqnarray}
which can be considered as a quantum generalization of the characteristic
polynomial (\ref{eq:char-pol-cl-3}). The ordering is important in
(\ref{eq:def-D-op}) since we are going to replace the parameters $u$, $z$,
${\cal H}_j$ by non-commuting operators.

\begin{prop}
 The operator $M$ (\ref{eq:def-M}) satisfies the equations
\beq
 XM-MH_3=0,
\label{eq:X-eq-op-M}
\eeq
\beq
 \D(y_j,Y_j;MH_1,MH_2,MH_3)=0, \qquad j=1,2
\label{eq:char-eq-M-op}
\eeq
\label{pr:char-eq-op-M}
where $H_{1,2,3}$ are the quantum Hamiltonians (\ref{eq:def-H-3}).
\end{prop}

{\bf Proof.} Though the equality (\ref{eq:X-eq-op-M}) is easy to derive
from the fact that
$M$ respects the constraint $x=t_3$, or directly from (\ref{eq:def-M}),
we shall proceed, however, in a more methodical fashion allowing
to prove both (\ref{eq:X-eq-op-M}) and (\ref{eq:char-eq-M-op}) in the
same way. Let us rewrite first the operator identities 
(\ref{eq:X-eq-op-M}) and (\ref{eq:char-eq-M-op}) for $M$ as algebraic
identities for the kernel $\M$ (\ref{eq:def-ker-M}).

We define the Lagrange adjoint Hamiltonians $H_k^*$ as
\begin{subeqnarray}
 H_1^*&=&T_1^{-1}v_{12}v_{13}+T_2^{-1}v_{21}v_{23}+T_3^{-1}v_{31}v_{32}, \\
 H_2^*&=&T_1^{-1}T_2^{-1}v_{13}v_{23}+T_1^{-1}T_3^{-1}v_{12}v_{32}+
T_2^{-1}T_3^{-1}v_{21}v_{31}, \\
 H_3^*&=&T_1^{-1}T_2^{-1}T_3^{-1},
\label{eq:def-H*}
\end{subeqnarray}

\beq
 \frac{1}{2\pi i}\oint\frac{dt}{t}\,f(t)(Hg)(t)=
 \frac{1}{2\pi i}\oint\frac{dt}{t}\,(H^*f)(t)g(t).
\label{eq:Lagrange-adj}
\eeq

In particular, $T_j^*=T_j^{-1}$. Considering $M$ in (\ref{eq:X-eq-op-M})
and (\ref{eq:char-eq-M-op}) as integral operator, we can use integration 
by parts and switch $H_k$ to the kernel $\M$ replacing them by $H_k^*$
according to (\ref{eq:Lagrange-adj}) which results in the $q$-difference
equations for $\M$:
\beq
 (X-H_3^*)\M=0,
\label{eq:X-eq-ker-M}
\eeq
\beq
 \D(y_j,Y_j;H_1^*,H_2^*,H_3^*)\M=0, \quad j=1,2.
\label{eq:char-eq-ker-M}
\eeq

While (\ref{eq:X-eq-ker-M}) is obvious, (\ref{eq:char-eq-ker-M}) needs
more consideration. Note that,
by virtue of (\ref{eq:act-Y-ypm}) and
(\ref{eq:act-T-tpm}), the action of ${\cal D}$ on $\M(y_+,y_-\mid t_-)$
is well defined. Note also
that the equations (\ref{eq:X-eq-ker-M}) and (\ref{eq:char-eq-ker-M})
are the quantum counterparts, resp., of the classical separated equations
(\ref{sep-eq-X}) and (\ref{eq:char-eq-cl}).

The next step is to notice that the kernel $\M$
 (\ref{eq:def-ker-M}) satisfies
the four first order $q$-difference equations 
\beq
 Y_jT_k\M=\aa_k(y_j)\M, \qquad j,k\in\{1,2\},
\label{eq:4-eq-M}
\eeq
where (compare to classical (\ref{eq:def-alpha-cl}))
\beq
 \aa_k(y)=
\frac{(1-q\ell^3y)(t_k-\ell t_3)(\ell t_3 y-t_{3-k})(qt_k-t_{3-k})}%
{\ell(1-y)(q\ell t_k-t_3)(q\ell^2t_3y-t_{3-k})(t_k-t_{3-k})},
\qquad k=1,2
\label{eq:def-alpha-q}
\eeq
which are verified directly from (\ref{eq:def-ker-M}) using the relations
(\ref{eq:shift-qx}). Note that (\ref{eq:4-eq-M}) is
the quantum counterpart of (\ref{eq:def-Y})--(\ref{eq:four-eq-Yy}).

{\bf Remark.}
It is easy to verify that the system (\ref{eq:4-eq-M}) is holonomic,
that is the operators $\aa_k(y_j)^{-1}Y_jT_k$ commute, provided
$y_j$ and $t_k$ are bound by (\ref{eq:cl-constraint}).

We proceed now to derive the third-order $q$-difference relations in $y_j$
(\ref{eq:char-eq-ker-M}) for $\M$ from the first-order relations 
(\ref{eq:4-eq-M}). The proof parallels that of theorem \ref{th:cl-char-eq}
for the classical case. Let us write down the equations 
(\ref{eq:char-eq-ker-M}) explicitely
\begin{eqnarray}
\lefteqn{
\kern-2mm\left[(1-qy_j)(1-q^2y_j)\ell^3Y_j^3
-(1-qy_j)(1-q^2\ell^2y_j)\ell^2Y_j^2H^*_1\right.} \nonumber \\
&&\kern-2mm\left.+(1-q\ell y_j)(1-q^2\ell^3y_j)\ell Y_j H^*_2
-(1-q\ell^3y_j)(1-q^2\ell^3y_j)H^*_3\right]\M=0,
\label{eq:cubic-eq-M}
\end{eqnarray}
then substitute into (\ref{eq:cubic-eq-M}) the expressions (\ref{eq:def-H*})
for $H_j^*$ and split the left-hand-side of (\ref{eq:cubic-eq-M}) into
two terms
\beq
 T_3^{-1}\check Z_1+Y_j\check Z_2=0
\label{eq:Z1+Z2_q}
\eeq
where
\begin{subeqnarray}
\lefteqn{\kern-2cm
 \check Z_1=\left[-(1-qy_j)(1-q^2\ell^2y_j)\ell^2Y_j^2v_{31}v_{32} 
-(1-q\ell^3y_j)(1-q^2\ell^3y_j)T_1^{-1}T_2^{-1}
\right.} \nonumber \\
&& \phantom{\left[\right.}\left. 
+(1-q\ell y_j)(1-q^2\ell^3y_j)\ell Y_j
(T_1^{-1}v_{12}v_{32}+T_2^{-1}v_{21}v_{31})
\right]\M, \\
\lefteqn{\kern-2cm
\check Z_2= \left[(1-y_j)(1-qy_j)\ell^3Y_j^2
+(1-\ell y_j)(1-q\ell^3y_j)\ell T_1^{-1}T_2^{-1}v_{13}v_{23}
\right.}\nonumber \\
&& \phantom{\left[\right.}\left. 
-(1-y_j)(1-q\ell^2y_j)\ell^2Y_j(T_1^{-1}v_{12}v_{13}+T_2^{-1}v_{21}v_{23})
\right]\M.
\label{eq:def-Zq}
\end{subeqnarray}

Introducing the notation
\beq
 \aa_{12}(y)\equiv\left.\aa_1(qy)\right|_{t_2:=qt_2}\!\aa_2(y)=
\left.\aa_2(qy)\right|_{t_1:=qt_1}\!\aa_1(y)\,,
\label{eq:def-alpha-12}
\eeq
\beq
 \w_{jk}=\frac{\ell^{-\half}t_j-q\ell^{\half}t_k}{t_j-qt_k}
\label{eq:def-wq}
\eeq
and noting that
\beq
  T_kv_{jk}=\w_{jk}T_k, \qquad  v_{jk}T_j=T_j\w_{jk},
\label{eq:ident-wT}
\eeq
it is easy to verify the algebraic identities for $\aa_{1,2}$
\begin{subeqnarray}
\lefteqn{\kern-2cm
-(1-qy)(1-q^2\ell^2y)\ell^2\w_{31}\w_{32}\aa_{12}(y)
-(1-q\ell^3y)(1-q^2\ell^3y)
} \nonumber \\
&& 
+(1-q\ell y)(1-q^2\ell^3y)\ell
\bigl(\w_{12}\w_{32}\aa_2(y)+\w_{21}\w_{31}\aa_1(y)\bigr)=0,
\\
\lefteqn{\kern-2cm
(1-y)(1-qy)\ell^3\aa_{12}(y)
+(1-\ell y)(1-q\ell^3y)\ell v_{13}v_{23}
} \nonumber \\
&& 
-(1-y)(1-q\ell^2y)\ell^2
\bigl(\w_{12}v_{13}\aa_2(y)+\w_{21}v_{23}\aa_1(y) \bigr)
=0.
\label{eq:alpha-id}
\end{subeqnarray}

Now we insert $T_kT_k^{-1}$ in appropriate places in (\ref{eq:def-Zq})
in such a way that $T_1^{-1}T_2^{-1}$ could be carried out to the left
of $[\ldots]$. Then we push the products $YT$ to the right using
(\ref{eq:ident-wT}) until they hit $\M$, so that (\ref{eq:4-eq-M}) could be
applied.
The equalities $\check Z_1=0$, $\check Z_2=0$
and therefore (\ref{eq:cubic-eq-M}) and (\ref{eq:char-eq-M-op})
follow then immediately from (\ref{eq:alpha-id}).
\endproof

\begin{prop}
The function $\bigl(MP_{\l}^{(\ell;q)}\bigr)(x,y_1,y_2)$ 
satisfies the $q$-difference equations (separated equations)
\beq
 (X-h_3)MP_{\l}^{(\ell;q)}=0,
\label{eq:XMJ}
\eeq
\beq
 \D(y_j,Y_j;h_1,h_2,h_3)MP_{\l}^{(\ell;q)}=0, \qquad j=1,2.
\label{eq:sep-eq-MJ}
\eeq
\label{pr:diff-eq-MJ}
\end{prop}
{\bf Proof.}
Apply the operator expressions (\ref{eq:X-eq-op-M}) and 
(\ref{eq:char-eq-M-op}) to the function 
$P_{\l}^{(\ell;q)}$. Using 
the operator ordering convention and the fact that Macdonald
polynomials $P_{\l}^{(\ell;q)}$ are the eigenfunctions of the 
Hamiltonians $H_j$ (\ref{eq:eigv-S}) one replaces $H_j$ by $h_j$.
Since $h_j$ are just numbers, the operator $M$ can be applied then 
directly to $P_{\l}^{(\ell;q)}$ which results in 
(\ref{eq:XMJ}) and (\ref{eq:sep-eq-MJ}). 
\endproof

In order to derive the factorization (\ref{eq:act-M-Macd-3}) of
$MP_{\l}^{(\ell;q)}$ we need more specific information about how
$M$ acts on the symmetric polynomials from ${\rm Sym}(t_1,t_2,t_3)$. 
Note that solutions to (\ref{eq:4-eq-M}), as to any $q$-difference
equations,
are defined only up to a factor invariant under $q$-shifts (quasiconstant).
Our choice (\ref{eq:def-ker-M}) of the kernel $\M$ corresponds to a
particular choice of the quasiconstant which is crucial for the results
given below.

Since the kernel $\M$ (\ref{eq:def-ker-M}) is a particular case
(\ref{eq:ab-g}) of the kernel $\M_{\a\b}$ (\ref{eq:def-ker-Mab}),
we can make use of the results obtained for $\M_{\a\b}$ in Appendix B.

Let us define few polynomial spaces. Let ${\rm Sym}(t_1,t_2,t_3)$ be
the space of Laurent polynomials symmetric w.r.t.\ permutations of 3 
variables $t_1$, $t_2$, $t_3$. A basis in ${\rm Sym}(t_1,t_2,t_3)$ is given
by $m_{\l}$ or $P_{\l}^{(\ell;q)}$. Let ${\rm Sym}(t_1,t_2;t_3)$
be the space of Laurent polynomials of the same 3 variables, symmetric only
w.r.t.\ $t_1\leftrightarrow t_2$. Obviously, 
${\rm Sym}(t_1,t_2;t_3)\supset{\rm Sym}(t_1,t_2,t_3)$.
Though the Macdonald polynomials belong to ${\rm Sym}(t_1,t_2,t_3)$
it is convenient to define $M$ on a larger space ${\rm Sym}(t_1,t_2;t_3)$.

Let ${\rm Ref}(t_-;t_+;t_3)$ be the space of Laurent polynomials in
$t_\pm$, $t_3$ which are reflexive in $t_-$ (invariant w.r.t.
$t_-\rightarrow t_-^{-1}$) and even in $t_\pm$ (invariant w.r.t.
$(t_-,t_+)\rightarrow(-t_-,\allowbreak -t_+)$.
Note that the change of variables
$(t_1,t_2,t_3)\rightarrow(t_-,t_+,t_3)$, see (\ref{eq:def-tpm}a), provides
an isomorphism ${\rm Sym}(t_1,t_2;t_3)\simeq{\rm Ref}(t_-;t_+;t_3)$.

The spaces ${\rm Sym}(y_1,y_2;x)\simeq{\rm Ref}(y_-;y_+;x)$
are defined similarly.

\begin{prop}
$$ M:{\rm Sym}(t_1,t_2;t_3)\rightarrow{\rm Sym}(y_1,y_2;x).$$

In particular, the image of a Macdonald polynomial
$P_{\l}^{(\ell;q)}\in{\rm Sym}(t_1,t_2,t_3)$ also lies in
${\rm Sym}(y_1,y_2;x)$.
\label{pr:act-M-poly}
\end{prop}

{\bf Proof.} The proposition \ref{pr:act-Mab-poly} from Appendix B
implies that $M:{\rm Ref}(t_-;t_+;t_3)\rightarrow{\rm Ref}(y_-;y_+;x)$.
Using the isomorphisms ${\rm Sym}(t_1,t_2;t_3)\simeq{\rm Ref}(t_-;t_+;t_3)$
and ${\rm Sym}(y_1,\allowbreak y_2;x)\simeq{\rm Ref}(y_-;y_+;x)$
we conclude the proof.
\endproof

Now everything is ready to prove the main statement of the theorem
\ref{th:main-th}.

\begin{prop}
 The operator $M$ transforms any Macdonald polynomial $P_{\l}^{(\ell;q)}$
into the product (\ref{eq:act-M-Macd-3}).
\label{pr:factorization-MJ}
\end{prop}
{\bf Proof.} We have already established that $MP_\l^{(\ell;q)}$
is a Laurent polynomial (proposition \ref{pr:act-M-poly}) satisfying
the $q$-difference equations (\ref{eq:XMJ}) and (\ref{eq:sep-eq-MJ}).
The factorization (\ref{eq:act-M-Macd-3}) follows from the fact
that $x^{h_3}$ and $S_{\l}^{(\ell;q)}(y)$ are the unique, up to
a constant factor, Laurent-polynomial solutions, of the $q$-difference
equations, resp.\ $(X-h_3)f(x)=0$, and $\D(y,Y;h_1,h_2,h_3)f(y)=0$.
The first statement is obvious, as for the second one,  see
proposition \ref{pr:unique-poly}.
\endproof

Though for the theorem \ref{th:main-th} we have used only the polynomiality
of $MP_\l^{(\ell;q)}$, in fact, the action of $M$ on 
${\rm Sym}(t_1,t_2;t_3)$ can be described in much more detail. Namely,
taking the formula (\ref{eq:act-Mab-pm}) from Appendix B,
making the substitutions (\ref{eq:ab-g}) and performing the changes of
variables (\ref{eq:def-tpm}a) and (\ref{eq:def-ypm}a) one
obtains the following result.

\begin{theo}
Consider the basis in ${\rm Sym}(t_1,t_2;t_3)$
\beq
p_{jk\nu}:=t_3^{j-2k}
t_1^kt_2^k(\ell^{-1}t_1t_3^{-1},\ell^{-1}t_2t_3^{-1};q)_\nu,
\qquad j,k\in\Z, \quad \nu\in\Z_{\geq0},
\label{eq:def-basis-t}
\eeq
and in ${\rm Sym}(y_1,y_2;x)$
\beq
 \tilde{p}_{jk\nu}:=x^jy_1^ky_2^k(y_1,y_2;q)_\nu,
\qquad j,k\in\Z, \quad \nu\in\Z_{\geq0},
\label{eq:def-basis-y}
\eeq
respectively.  The operator $M$ acts on $p_{jk\nu}$ as follows
\beq
 M:p_{jk\nu}\rightarrow \ell^{3k}\frac{(\ell^{-2};q)_\nu}{(\ell^{-3};q)_\nu}
\tilde{p}_{jk\nu}.
\label{eq:act-M-basis-t}
\eeq
\label{th:act-M-basis-t}
\end{theo}

Postponing the proof of the formulas (\ref{eq:def-sep-poly-3})
and (\ref{eq:def-chi-3}) for the next section, we can prove now the final
statement of theorem \ref{th:main-th}.

\begin{prop}
The normalization coefficient $c_{\l}$ in (\ref{eq:act-M-Macd-3})
is given by (\ref{eq:def-c-3}).
\label{pr:expr-cn}
\end{prop}
{\bf Proof.}
In this case, it is convenient to make use of the isomorphisms described
above and to think of $M$ as acting from ${\rm Ref}(t_-;t_+;t_3)$ into
${\rm Ref}(y_-;y_+;x)$. 
Comparing the asymptotics of
the monomial symmetric functions $m_\l$
\be
 m_{\l_1\l_2\l_3}\sim t_-^{\l_3-\l_1}t_+^{\l_3+\l_1}t_3^{\l_1+\l_2+\l_3},
\qquad t_-\rightarrow\infty
\ee
and of the polynomial $p_{jk\nu}$ (\ref{eq:def-basis-t})
\begin{eqnarray*}
 p_{jk\nu}&\equiv& t_3^{j}t_+^{2k}
(\ell^{-1}t_+t_-,\ell^{-1}t_+t_-^{-1})_\nu \\
&\sim& (-1)^\nu q^{\frac{\nu(\nu-1)}{2}}\ell^{-\nu}
t_3^{j}t_+^{2k+\nu}t_-^\nu, \qquad t_-\rightarrow\infty
\end{eqnarray*}
we conclude that the transition matrix between the bases $m_\l$ and
$p_{jk\nu}$ is triangular
\beq
 m_\l=(-1)^{\l_{31}}q^{-\frac{\l_{31}(\l_{31}-1)}{2}}\ell^{\l_{31}}
p_{|\l|,\l_1,\l_{31}}
+\sum_{\nu<\l_{31}}\sum_{j,k}(\ldots)p_{jk\nu}.
\label{eq:triang-m-p}
\eeq

Given the mutual triangularity (\ref{eq:def-Macd-poly}) of the bases
$P_\l^{(\ell;q)}$ and $m_\l$, it means that the expansion of $P_\l^{(\ell;q)}$
in $p_{jk\nu}$ has the same structure as (\ref{eq:triang-m-p}).
Using then (\ref{eq:act-M-basis-t})
and the asymptotics of $\tilde{p}_{jk\nu}$
(\ref{eq:def-basis-y})
\be
 \tilde{p}_{jk\nu}\equiv x^j y_+^{2k}(y_+y_-,y_+y_-^{-1})_\nu 
\sim (-1)^\nu q^{\frac{\nu(\nu-1)}{2}}x^jy_+^{2k+\nu}y_-^\nu,
\qquad y_-\rightarrow\infty
\ee
we obtain
\begin{eqnarray}
 MP_\l^{(\ell;q)}&=&
(-1)^{\l_{31}}q^{-\frac{\l_{31}(\l_{31}-1)}{2}}\ell^{2\l_1+\l_3}
\frac{(\ell^{-2};q)_{\l_{31}}}{(\ell^{-3};q)_{\l_{31}}}
\tilde{p}_{\abs{\l},\l_1,\l_{31}} +\ldots \nonumber \\
&\sim&\ell^{2\l_1+\l_3}
\frac{(\ell^{-2};q)_{\l_{31}}}{(\ell^{-3};q)_{\l_{31}}}
  x^{\abs{\l}}y_+^{\l_3+\l_1}y_-^{\l_{31}}, \qquad y_-\rightarrow\infty.
\label{eq:asympt-MP}
\end{eqnarray}

On the other hand, (\ref{eq:pn-poly}) implies that
\be
 c_{\l}\;x^{\abs{\l}}\;
S_{\l}^{(\ell;q)}(y_+y_-)\;S_{\l}^{(\ell;q)}(y_+y_-^{-1})
 \sim c_{\l}\;\chi_{\l,\l_1}^{(\ell;q)}\;\chi_{\l,\l_3}^{(\ell;q)}\;
x^{\abs{\l}}y_+^{\l_3+\l_1}y_-^{\l_3-\l_1}
\ee
whence
\beq
 c_{\l}\;\chi_{\l,\l_1}^{(\ell;q)}\;\chi_{\l,\l_3}^{(\ell;q)}=
\frac{(\ell^{-2};q)_{\l_3-\l_1}}{(\ell^{-3};q)_{\l_3-\l_1}}\ell^{\l_3+2\l_1}.
\label{eq:def-cn}
\eeq

It remains only to use the formulas (\ref{eq:chi-n1-n3})
proved in the next section, and obtain (\ref{eq:def-c-3}).
\endproof

Compared to \cite{Kuz-Skl} our formula (\ref{eq:def-c-3}) for the 
normalization coefficients $c_\l$ is new, and its nonrelativistic analog
\be
   c_\l=\frac{(2g)_{\l_{31}}(2g)_{\l_{32}}(g)_{\l_{21}}}%
{(3g)_{\l_{31}}(g)_{\l_{32}}(2g)_{\l_{21}}}\,,\qquad
(\a)_k\equiv \a(\a+1)\ldots(\a+k-1)\,,
\ee
fills the gap in the description given in \cite{Kuz-Skl}
of the integral representation for Jack polynomials analogous to
(\ref{eq:int-repr-macd}).

We conclude this section with a list of results concerning the inverse
operator $M^{-1}$. All the preparatory work being done in Appendix B,
it remains only to use the correspondence (\ref{eq:ab-g}) between
$M_{\a\b}$ and $M$.

{}From (\ref{eq:inv-Mab}) and (\ref{eq:def-ker-Mab-inv}) it follows that
$M^{-1}$ is an integral operator
\begin{eqnarray}
\lefteqn{
 M^{-1}: \Phi(x,y_1,y_2)\rightarrow\Psi(t_1,t_2,t_3) }\nonumber \\
&&\kern-12mm=\frac{1}{2\pi i}
\int\limits_{\G_{-g,3g}^{t_+,t_-}} \!\!\frac{dy_-}{y_-}
\tilde\M\!\left(\frac{(t_1t_2)^{\frac12}}{t_3},
             \left(\frac{t_1}{t_2}\right)^{\frac12}\Bigm| y_-\right)
\!\Phi\!\left(t_3,\frac{\ell^{-\frac32}(t_1t_2)^{\frac12}y_-}{t_3},
     \frac{\ell^{-\frac32}(t_1t_2)^{\frac12}}{t_3y_-}\right)
\label{eq:def-M1}
\end{eqnarray}
with the kernel 
\beq
\tilde\M(t_+,t_-\mid y_-)=\frac{(1-q)(q;q)_\infty^2(y_-^2,y_-^{-2};q)_\infty
\,\Lq{\ell^{-1}}{t_-}{t_+} }%
{2B_q(-g,3g)
\,\Lq{\ell^{\half}}{y_-}{t_-}\Lq{\ell^{-\frac32}}{y_-}{t_+}  }.
\label{eq:def-ker-M-inv}
\eeq

Reversing (\ref{eq:act-M-basis-t}) one obtains the formula for the action of
$M^{-1}$ on the basis $\tilde{p}_{jk\nu}$
\beq
 M^{-1}: \tilde{p}_{jk\nu}\rightarrow 
\ell^{-3k}\frac{(\ell^{-3};q)_\nu}{(\ell^{-2};q)_\nu} p_{jk\nu}.
\label{act-M1-basis-y}
\eeq

Reversing (\ref{eq:act-M-Macd-3}) provides a new integral representation
of $A_2$ Macdonald polynomials in terms of the Laurent polynomials
$S_\l^{(\ell;q)}(y)$ (\ref{eq:def-sep-poly-3})
\beq
 M^{-1}: c_{\l}x^{\abs{\l}}S_{\l}^{(\ell;q)}(y_1)S_{\l}^{(\ell;q)}(y_2)
\rightarrow P_{\l}^{(\ell;q)}(t_1, t_2, t_3).
\label{eq:int-repr-macd}
\eeq

Finally, from the propositions \ref{pr:Mab-int}  and
\ref{pr:inv-Mab} it follows that for 
positive integer $g$ the operator $M^{-1}$ turns into a $q$-difference
operator of order $g$:
\beq
 M^{-1}:\Phi(x,y_1,y_2)\rightarrow\sum_{k=1}^g
\xi_k\!\left(\frac{(t_1t_2)^{\half}}{t_3},
 \left(\frac{t_1}{t_2}\right)^{\half}\right)
\Phi\!\left(t_3,q^{g+k}\frac{t_1}{t_3},q^{2g-k}\frac{t_2}{t_3}\right)
\label{eq:M1-int-g}
\eeq
where $\xi_k(r,s)$ is given by (\ref{eq:expr-xik}).
The result is not surprising in view of the similar result for the
nonrelativistic case \cite{Kuz-Skl} where $M^{-1}$ becomes a differential
operator of order $g$ for $g\in\Z_{\geq0}$. In \cite{Kuz-Skl}
this result was derived using a representation of $M^{-1}$ in terms of the
fractional differentiation operator. In the relativistic case it is
also possible to relate $M^{-1}$ with a sort of fractional $q$-difference
operator. We intend to touch this subject in a separate paper.

\newsection{Separated equation}

In this section the results are collected concerning the Laurent
polynomials $S_{\l}^{(\ell;q)}(y)$ and the corresponding 
$q$-difference equations. Since all the results are easy to generalize 
from $n=3$ to  arbitrary $n$, we give them in the most general form.
\begin{guess}
The correct generalization of the formula (\ref{eq:def-sep-poly-3})
for $S^{(\ell;q)}_{\l}(y)$ for any $n$ is given by 
\beq
 S^{(\ell;q)}_{\l}(y)=y^{\l_1}(y;q)_{1-ng}\,
   {}_n\/\oldphi_{n-1}\left[\begin{array}{c}
     a_1,\ldots,a_n \\
     b_1,\ldots,b_{n-1} \end{array};q,y\right]
\label{eq:def-sep-poly-gen}
\eeq
where
\beq
a_j=\ell^{n-j+1}q^{\l_1-\l_{n-j+1}+1}, \qquad b_j=a_j\ell^{-1}.
\label{eq:def-aj-bj}
\eeq
\end{guess}

\begin{prop}
$S^{(\ell;q)}_{\l}(y)$ is a Laurent polynomial in $y$ of the form
\beq
 S^{(\ell;q)}_{\l}(y)=\sum_{k=\l_1}^{\l_n} \chi^{(\ell;q)}_{\l,k}y^k.
\label{eq:def-Sn}
\eeq
\end{prop}
{\bf Proof.} Observe, first, that if 
$a=bq^\nu$ for some positive integer $\nu$ then
\beq
 \frac{(a;q)_k}{(b;q)_k}=\frac{(bq^k;q)_\nu}{(b;q)_\nu}
\label{eq:lemma-ab-1}
\eeq
is a polynomial in $q^k$ of degree $\nu$ whose coefficients are rational
functions in $b$ and $q$.
As a consequence, if $a_{j+1}=b_jq^{\nu_j}$ then
\beq
 P_N(q^k)\equiv
 \frac{(a_2;q)_k\ldots(a_n;q)_k}{(b_1;q)_k\ldots(b_{n-1};q)_k}=
\frac{(b_1q^k;q)_{\nu_1}\ldots(b_{n-1}q^k;q)_{\nu_{n-1}}}%
{(b_1;q)_{\nu_1}\ldots(b_{n-1};q)_{\nu_{n-1}}}
\label{eq:lemma-ab-2}
\eeq
is a polynomial in $q^k$ of degree $N=\nu_1+\cdots+\nu_{n-1}$.

In our case, $\nu_j=\l_{n-j+1}-\l_{n-j}$, $N=\l_n-\l_1$ by virtue of 
(\ref{eq:def-aj-bj}), and from (\ref{eq:def-sep-poly-gen}) and 
(\ref{eq:def-bas-hgs}) one obtains
\beq
 {}_n\/\oldphi_{n-1}\left[\begin{array}{c}
     a_1,\ldots,a_n \\
     b_1,\ldots,b_{n-1} \end{array};q,y\right]=
 \sum_{k=0}^\infty\frac{(a_1;q)_k}{(q;q)_k}\,y^k P_N(q^k)
\label{eq:def-PM}
\eeq
where $P_N(q^k)$ is given by (\ref{eq:lemma-ab-2}). It remains now to apply 
the following lemma.

\begin{lemma}
Let $P_N(y)$ be a polynomial in $y$ of degree $\leq N$. Then
\beq
 \sum_{k=0}^\infty\frac{(a;q)_k}{(q;q)_k}\,y^k P_N(q^k)=
Q_N(y)\frac{(aq^Ny;q)_\infty}{(y;q)_\infty}
\label{eq:PQ}
\eeq
where $Q_N(y)$ is a polynomial in $y$ of degree $\leq N$.
\label{lemma:PQ}
\end{lemma}

{\bf Proof.} It is sufficient to consider the polynomials 
$P_N(q^k)=(q^{k-\nu+1};q)_\nu$ for $\nu=0,1,\ldots,N$
forming a basis in the polynomial ring.
Then
\begin{eqnarray}\lefteqn{
 \sum_{k=0}^\infty\frac{(a;q)_k}{(q;q)_k}\,y^k (q^{k-\nu+1};q)_\nu=
 \sum_{k=\nu}^\infty[\cdots]= 
 \sum_{k=\nu}^\infty\frac{(a;q)_k}{(q;q)_{k-\nu}}\,y^k} \nonumber \\
&&  =\sum_{k=0}^\infty\frac{(a;q)_{k+\nu}}{(q;q)_k}\,y^{k+\nu}=
  (a;q)_\nu y^\nu  \sum_{k=0}^\infty\frac{(aq^\nu;q)_k}{(q;q)_{k}}\,y^k. 
\label{eq:lemma-2.2}
\end{eqnarray}

Using the formula (\ref{eq:1phi0}) and the identity 
$(aq^\nu;q)_\infty=(aq^\nu;q)_{N-\nu}(aq^N;q)_\infty$
one obtains finally the expression (\ref{eq:PQ}) where
$ Q_N(y)=(a;q)_\nu y^\nu(aq^\nu y;q)_{N-\nu}$.
\endproof

Applying the above lemma to the case of the polynomial $P_N(q^k)$ given by 
(\ref{eq:def-PM}) and $a=a_1=\ell^nq^{\l_1-\l_n+1}$
we obtain finally that $y^{-\l_1}S^{(\ell;q)}_{\l}(y)$ is a polynomial 
of degree 
$\leq \l_n-\l_1$. 
\endproof

\begin{prop}
The coefficients $\chi^{(\ell;q)}_{\l,k}$ in the expansion 
(\ref{eq:def-Sn}) are given by
\beq
\chi_{\l,k}^{(\ell;q)}=\left(q\ell^n\right)^{k-\l_1}
\frac{(q^{-1}\ell^{-n};q)_{k-\l_1}}{(q;q)_{k-\l_1}}\,
{}_{n+1}\oldphi{}_{n}\!\left[\tworow{q^{\l_1-k},a_1,\ldots,a_n}%
{q^{\l_1-k+2}\ell^n,b_1,\ldots,b_{n-1}};q,q\right].
\label{eq:chi-N}
\eeq

In particular, for $n=3$, (\ref{eq:chi-N}) produces (\ref{eq:def-chi-3}).
\label{pr:chi-N}
\end{prop}
{\bf Proof.} We know already that $S^{(\ell;q)}_{\l}(y)$ is a Laurent
polynomial and is thus defined for any $y\in\C\setminus\{0,\infty\}$.
Suppose for a while that $\abs{y}<1$. Then both factors $(y;q)_{1-ng}$
and ${}_n\oldphi_{n-1}$ in (\ref{eq:def-sep-poly-gen}) are given by the
 convergent series  (\ref{eq:1phi0}) and
(\ref{eq:def-bas-hgs}), respectively. Multiplying the two power series
in $y$ we observe that the coefficients at $y^k$ is expressed in terms
of ${}_{n+1}\oldphi_n$ series:
\beq
 S_{\l}^{(\ell;q)}(y)=y^{\l_1}\sum_{k=0}^{\infty}
\left(q\ell^n\right)^{k}
\frac{(q^{-1}\ell^{-N};q)_{k}}{(q;q)_{k}}
\;{}_{n+1}\oldphi{}_{n}\left[\tworow{q^{-k},a_1,\ldots,a_n}%
{q^{-k+2}\ell^n,b_1,\ldots,b_{n-1}};q,q\right]\;y^k.
\label{eq:S-inf-ser}
\eeq

In fact, the sum in (\ref{eq:S-inf-ser}) is finite: $\sum_{k=0}^{\l_{n1}}$.
To see this, use the formula (1.9.11) from \cite{GR}: let $\nu$, $k_1,\ldots,
k_n\in\Z_{\geq0}$, then
\beq
{}_{n+1}\oldphi{}_n\left[\tworow{q^{-\nu},b_1q^{k_1},\ldots,b_nq^{k_n}}%
{b_1,\ldots,b_n};q,q\right]=0
\eeq
for $\nu>k_1+\cdots+k_n$. Substituting
\be
\nu=k\,, \qquad b_n=q^{2-k}\ell^n\,,\qquad
k_n=\l_1-\l_n+k-1\,,
\ee
\be
b_j=\ell^{n-j}q^{\l_1-\l_{n-j+1}+1}\,,\qquad 
k_j=\l_{n-j+1}-\l_{n-j}, \qquad j=1,\ldots,n-1\,,
\ee
we obtain that
\be
{}_{n+1}\oldphi{}_n\left[\tworow{q^{-k},a_1,\ldots,a_n}%
{q^{2-k}\ell^n,b_1,\ldots,b_{n-1}};q,q\right]=0
\ee
for $k\geq \l_n-\l_1+1$, hence 
the sum in (\ref{eq:S-inf-ser}) is finite: $\sum_{k=0}^{\l_{n1}}$. The
coefficient at $y^k$ in (\ref{eq:S-inf-ser}) produces, respectively, 
(\ref{eq:chi-N}).
\endproof

For the sake of reference we present a short list of polynomials
$S_\l^{(\ell;q)}(y)$ in case $n=3$:
\be
S_{000}^{(\ell;q)}=1,\quad
S_{001}^{(\ell;q)}=1+{\tfrac {{\ell}^{2}y}{\ell+1}},\quad
S_{011}^{(\ell;q)}=1+\ell\,\left (\ell+1\right )y,
\ee
\be
S_{002}^{(\ell;q)}=1
+\tfrac{\ell^2(q\ell+\ell-q-1)y}{\ell^2-q}
+\tfrac{(\ell-q)\ell^4y^2}{(\ell^2-q)(\ell+1)}\,,
\ee
\be
S_{012}^{(\ell;q)}=1
+\tfrac{(\ell^3+\ell^2q+\ell^2-\ell q-\ell-q)\ell y}{\ell^2-q}
+\ell^3y^2\,,
\ee
\be
S_{022}^{(\ell;q)}=1
+\tfrac{(\ell^2-1)(q+1)\ell y}{\ell-q}
+\tfrac{(\ell^2-q)(\ell+1)\ell^2y^2}{\ell-q}\,,
\ee
\be
S_{003}^{(\ell;q)}=1
+\tfrac{(1+q+q^2)(\ell-1)\ell^2 y}{\ell^2-q^2}
+\tfrac{(1+q+q^2)(\ell-1)\ell^4 y^2}{(\ell+q)(\ell^2-q)}
+\tfrac{(\ell-q^2)\ell^6 y^3}{(\ell+1)(\ell+q)(\ell^2-q)}\,,
\ee
\be
S_{013}^{(\ell;q)}=1
+\tfrac{(\ell^3+q^2\ell^2+\ell^2q+\ell^2-q^2\ell-\ell q-\ell-q^2)\ell y}
{\ell^2-q^2}
\ee\be
+\tfrac{(q\ell^3+\ell^3+\ell^2q^2+\ell^2q+\ell^2
-q^3\ell-\ell q^2-\ell q-q^3-q^2)(\ell-1)\ell^3 y^2}
{(\ell^2-q)(\ell^2-q^2)}
+\tfrac{(\ell-q)\ell^5 y^3}{\ell^2-q}\,,
\ee
\be
S_{023}^{(\ell;q)}=1
+\tfrac{(q\ell^3+\ell^3+\ell^2q^2+\ell^2q+\ell^2-q^3\ell
-q^2\ell-\ell q-q^3-q^2)(\ell-1)\ell y}
{(\ell-q)^2(\ell+q)}
\ee\be
+\tfrac{(\ell^3+\ell^2q^2+\ell^2q+\ell^2
-q^2\ell-\ell q-\ell-q^2)(\ell^2-q)\ell^2 y^2}
{(\ell-q)^2(\ell+q)}
+\tfrac{(\ell^2-q)\ell^4 y^3}{\ell-q}\,,
\ee
\begin{eqnarray*}
S_{033}^{(\ell;q)}&\!=\!&1
+\tfrac{(1+q+q^2)(\ell^2-1)\ell y}{\ell-q^2}
+\tfrac{(1+q+q^2)(\ell^2-1)(\ell^2-q)\ell^2 y^2}{(\ell-q)(\ell-q^2)} \\
&&+\tfrac{(\ell+1)(\ell+q)(\ell^2-q)\ell^3 y^3}{\ell-q^2}\,. 
\end{eqnarray*}

{\bf Remark.} It easy to give more simple expressions for some of
$\chi_{\l,k}^{(\ell;q)}$ such as
\beq
 \chi_{\l,\l_1}^{(\ell;q)}=1,
\label{eq:chi-n1}
\eeq
and
\beq
\chi_{\l,\l_n}^{(\ell;q)}=\ell^{\abs{\l}-n\l_1}\prod_{j=1}^{n-1}
\frac{(\ell^{-j};q)_{\l_j-\l_1}(\ell^{-j};q)_{\l_n-\l_{n-j}}}%
{(\ell^{-j};q)_{\l_{j+1}-\l_1}(\ell^{-j};q)_{\l_n-\l_{n-j+1}}}.
\label{eq:chi-nN}
\eeq
which for $n=3$ produce
\beq
\chi_{\l,\l_1}^{(\ell;q)}=1, \qquad
\chi_{\l,\l_3}^{(\ell;q)}=\ell^{-2\l_1+\l_2+\l_3}
\frac{(\ell^{-1};q)_{\l_{32}}(\ell^{-2};q)_{\l_{21}}}%
{     (\ell^{-2};q)_{\l_{32}}(\ell^{-1};q)_{\l_{21}}}.
\label{eq:chi-n1-n3}
\eeq 

The formula (\ref{eq:chi-n1}) is obvious. To obtain (\ref{eq:chi-nN}), use
the summation formula (1.9.10) from \cite{GR}:
\begin{eqnarray}
\lefteqn{
{}_{n+1}\oldphi_{n}\left[
\tworow{q^{-\nu},\b,\b_1q^{k_1},\ldots,\b_{n-1}q^{k_{n-1}}}%
{\b q,\b_1,\ldots,\b_{n-1}};q,q
\right] }\nonumber \\
&&=
\frac{(q;q)_\nu(\b_1/\b;q)_{k_1}\ldots(\b_{n-1}/\b;q)_{k_{n-1}}}%
{(\b q;q)_\nu(\b_1;q)_{k_1}\ldots(\b_{n-1};q)_{k_{n-1}}}\,\b^\nu,
\end{eqnarray}
where
$\nu,k_1,\ldots,k_{n-1}\in\Z_{\geq0}$ and $\nu\geq k_1+\cdots+k_{n-1}$.
Substituting
\be
  \nu=\l_{n1},\qquad \b=\ell^n q^{1-\l_{n1}}\equiv a_1\,,
\ee
\be
\b_j=\ell^j q^{1-\l_{j+1}+\l_1}\equiv b_{n-j}, \qquad
k_j=\l_{j+1}-\l_j,\quad j=1,\ldots,n-1\,,
\ee
we obtain, after a series of equivalent transformations
(see Appendix I to \cite{GR}), the expression
(\ref{eq:chi-nN}).

{\bf Remark.} There is also a simple formula for 
$S^{(\ell;q)}_{\l}(\ell^{-n})$:
\beq
 S^{(\ell;q)}_{\l}(\ell^{-n})=\ell^{-n\l_1}(\ell^{-n};q)_{\l_{n1}}
\prod_{j=1}^{n-1}
\frac{(\ell^{-j};q)_{\l_j-\l_1}}{(\ell^{-j};q)_{\l_{j+1}-\l_1}},
\eeq 
or, for $n=3$,
\beq
S^{(\ell;q)}_{\l}(\ell^{-3})=\ell^{-3\l_1}
\frac{(\ell^{-2};q)_{\l_{21}}(\ell^{-3};q)_{\l_{31}}}%
{(\ell^{-1};q)_{\l_{21}}(\ell^{-2};q)_{\l_{31}}},
\eeq
which are proved in a way similar to (\ref{eq:chi-nN}) using the 
formula (1.9.9) from \cite{GR}.

The polynomials $S^{(\ell;q)}_{\l}(y)$ can be expressed  also 
in terms of the $q$-Lauricella function (\ref{eq:def-q-Lauri}).
\begin{prop}
The following equalities hold:
\begin{eqnarray}
S_{\l}^{(\ell;q)}(y)&=&y^{\l_1}\frac{(q\ell^nq^{\l_{1n}}y;q)_{\l_{n1}}}
{\prod_{j=1}^{n-1}(q^{\l_1-\l_{n-j+1}+1}\ell^{n-j};q)_{\l_{n-j+1}-\l_{n-j}}}
\nonumber \\
&&\times \oldphi_D\left[\begin{array}{l}
a';b'_1,\ldots,b'_{n-1}\cr c\end{array};q;x_1,\ldots,x_{n-1}
\right]\,,
\label{eq:repr-S-Lauri}
\end{eqnarray}
where $\l_{ij}\equiv \l_i-\l_j$ and
\beq
a'=y, \quad c=q\ell^nq^{\l_{1n}}y, \quad
x_j=q\ell^{n-j}q^{\l_1-\l_{n-j}}, \quad b'_j=q^{\l_{n-j}-\l_{n-j+1}},
\label{eq:param-Lauri}
\eeq
for $j=1,\ldots,n-1$. Another expression for $S_{\l}^{(\ell;q)}(y)$ reads
\beq
S_{\l}^{(\ell;q)}(y)=y^{\l_1}(q\ell y;q)_{(1-n)g}
\left(\prod_{i=1}^{n-1}(a_i;q)_g\right)
\;\oldphi_D\left[\begin{array}{l}
y;\ell^{-1},\ldots,\ell^{-1}\cr
q\ell y\end{array};q;a_1,\ldots,a_{n-1}\right]\,.
\label{eq:repr-S-Lauri-2}
\eeq
\label{pr:Lauricella}
\end{prop}
{\bf Proof.}
The formula (\ref{eq:repr-S-Lauri}) is obtained by substituting the
parameters (\ref{eq:param-Lauri}) into Andrews' formula 
(\ref{eq:Andrews-Lauri})
for $q$-Lauricella function and comparing the result to
(\ref{eq:def-sep-poly-gen}). Note that $c/a'\equiv a_1$, 
$x_j\equiv a_{j+1}$, $b_j'x_j\equiv b_j$.

Similarly, substituting into (\ref{eq:Andrews-Lauri}) the parameters
$a'=y$, $c=q\ell y$, $b_j'=\ell^{-1}$, $x_j=a_j$ ($j=1,\ldots,n-1$)
such that $c/a'\equiv a_n$, $b_j'x_j\equiv b_j$, one arrives at 
(\ref{eq:repr-S-Lauri-2}).
\endproof

{\bf Corollary 1.}
Substituting the definition 
(\ref{eq:def-q-Lauri}) of $\oldphi_D$ into formula (\ref{eq:repr-S-Lauri})
we obtain another explicit representation for $S_{\l}^{(\ell;q)}(y)$:
\be
S_{\l}^{(\ell;q)}(y)=y^{\l_1}\frac{1}
{\prod_{j=1}^{n-1}(q^{\l_1-\l_{n-j+1}+1}\ell^{n-j};q)_{\l_{n-j+1}-\l_{n-j}}}
\ee
\be
\times \sum_{k_1=0}^{\l_n-\l_{n-1}}\;\cdots\; 
\sum_{k_{n-1}=0}^{\l_2-\l_{1}}\;(q\ell^nq^{\l_{1n}+k_1+\ldots+k_{n-1}}y;q)_{
\l_{n1}-k_1-\ldots-k_{n-1}}\;(y;q)_{k_1+\ldots+k_{n-1}}
\ee
\beq
\times \prod_{j=1}^{n-1} \frac{(q^{\l_{n-j}-\l_{n-j+1}};q)_{k_j}
(q\ell^{n-j}q^{\l_1-\l_{n-j}})^{k_j}}
{(q;q)_{k_j}}\,.
\label{eq:another-repr-S}
\eeq

{\bf Corollary 2.}
It is possible also to represent $S_{\l}^{(\ell;q)}(y)$
as a $q$-integral (\ref{eq:def-q-int}):
\be
S_{\l}^{(\ell;q)}(q^x)=q^{\l_1x}\frac{(q\ell^nq^{\l_{1n}}q^x;q)_{\l_{n1}}}
{\prod_{j=1}^{n-1}(q^{\l_1-\l_{n-j+1}+1}\ell^{n-j};q)_{\l_{n-j+1}-\l_{n-j}}}
\ee
\beq
\times \frac{1}{B_q(x,\l_{1n}+1-ng)}
\int_0^1 d_qt \;t^{x-1}\frac{(tq;q)_{\l_{1n}-ng}}
{\prod_{j=1}^{n-1}(tq\ell^{n-j}q^{\l_1-\l_{n-j}};q)_{\l_{n-j}-\l_{n-j+1}}}\,.
\label{eq:repr-S-q-int}
\eeq

To obtain the formula (\ref{eq:repr-S-q-int}) rewrite Andrews' formula
(\ref{eq:Andrews-Lauri}) as a $q$-integral
\beq
\oldphi_D\left[\begin{array}{l}
q^\alpha;q^{\beta_1},\ldots,q^{\beta_{n-1}} \cr
q^\gamma\end{array};q;x_1,\ldots,x_{n-1}\right]=
\frac{1}{B_q(\alpha,\gamma-\alpha)}
\int_0^1 d_qt \;t^{\alpha-1}\frac{(tq;q)_{\gamma-\alpha-1}}
{\prod_{j=1}^{n-1}(tx_j;q)_{\beta_j}}\,.
\eeq
and substitute
\be
\alpha=x\;\;\;(y=q^x), \quad \gamma=\l_{1n}+1-ng+x, \quad
\beta_j=\l_{n-j}-\l_{n-j+1}, \quad x_j=q\ell^{n-j}q^{\l_1-\l_{n-j}}\,.
\ee

The rest of the results are concerned with the separated $q$-difference
equations for the polynomials $S_{\l}^{(\ell;q)}(y)$.

\begin{prop}
The polynomial $f(y):=S^{(\ell;q)}_{\l}(y)$ (\ref{eq:def-sep-poly-gen}) 
satisfies the $q$-difference equation
\beq
 \sum_{k=0}^n(-1)^k\ell^{\frac{n-1}{2}k}(1-q^k\ell^ky)(y;q)_k\,
(q^{k+1}\ell^ny;q)_{n-k}\,h_{n-k}\,f(q^ky)=0
\label{eq:sep-eq-N}
\eeq
where, $h_k$ are given by (\ref{eq:def-hk-gen}) and,
as in the classical case (\ref{eq:char-poly-N}), 
we assume $h_0\equiv1$.
\label{theo:sep-eq-N}
\end{prop}

{\bf Proof.} Using the definitions (\ref{eq:def-sep-poly-gen}) and
(\ref{eq:def-aj-bj}) together with (\ref{eq:def-hk-gen}),
it is a matter of straightforward calculation to transform
the $q$-difference equation (\ref{eq:hg-diff-eq})
for ${}_n\oldphi_{n-1}$ into (\ref{eq:sep-eq-N}).
\endproof

In fact, the factor $(1-q^n\ell^ny)$ can be cancelled from 
(\ref{eq:sep-eq-N}) which results in the equation
\begin{eqnarray}
\lefteqn{
 (-1)^n\ell^{\frac{n(n-1)}{2}}(y;q)_nh_0f(q^ny) } \nonumber \\
&& \kern-6mm
+\sum_{k=0}^{n-1}(-1)^k\ell^{\frac{n-1}{2}k}(1-q^k\ell^ky)(y;q)_k\,
(q^{k+1}\ell^ny;q)_{n-k-1}\,h_{n-k}\,f(q^ky)=0\,.
\label{eq:sep-eq-N-simp}
\end{eqnarray}
In the case $n=3$ the $q$-difference equation (\ref{eq:sep-eq-N})
takes the form 
$$\D(y,Y;h_1,h_2,h_3)f(y)=0$$
 where $\D$ is given by (\ref{eq:def-D-op}), or, explicitely,
\begin{eqnarray}
\lefteqn{\kern-6mm (1-qy)(1-q^2y)\ell^3\,f(q^3y)
-(1-qy)(1-q^2\ell^2y)\ell^2h_1\,f(q^2y) 
}\nonumber \\
&&\kern-11mm +(1-q\ell y)(1-q^2\ell^3y)\ell h_2\,f(qy) 
-(1-q\ell^3y)(1-q^2\ell^3y)h_3\,f(y)=0. 
\label{eq:sep-eq-3}
\end{eqnarray}
\begin{prop}
Let 
\be
 G_n^{(0)}:=\Z\cup\half\Z\cup\ldots\cup\frac{1}{n-1}\Z, \qquad
 G_n^{(1)}:=\left\{\frac{1}{n},\frac{2}{n},\ldots,\frac{n-2}{n}\right\}.
\label{eq:def-G}
\ee

Then, for all $g>0$ except for the finite number of points
 $g\in G_n\equiv G_n^{(0)}\cap G_n^{(1)}$,
the separated polynomial $f(y):=S^{(\ell;q)}_{\l}(y)$ 
(\ref{eq:def-sep-poly-gen}) is the only, up to
a constant factor, Laurent-polynomial solution to the $q$-difference 
equation (\ref{eq:sep-eq-N}).

In particular, $G_3=\emptyset$, so for $n=3$ the uniqueness 
of L.-p.\ solution holds $\forall g>0$.
\label{pr:unique-poly}
\end{prop}
{\bf Proof.} 
In the nonrelativistic case \cite{Kuz-Skl} the analog of the equation
(\ref{eq:sep-eq-N}) is a differential equation having 3 regular 
singularities: ${0,1,\infty}$, and the uniqueness of L.-p.\ solution
is proved by analysis of the corresponding characteristic exponents.
As shown below, the argument can be translated rather directly to the
$q$-difference case.

Let $f(y)$ be a non-zero Laurent-polynomial solution to
(\ref{eq:sep-eq-N}) or, equivalently, (\ref{eq:sep-eq-N-simp}). 
Then, subsituting into (\ref{eq:sep-eq-N-simp}) the
values $y=q^{-j}$, $j=0,1,2,\ldots$, one observes that $f(q^{-j})$ can
be determined recursively, starting from $f(1)$ since the factor
$(y;q)_k$ cuts away the terms with $k>j$. The only obstacle could be
the vanishing of the factor $(q^{k+1}\ell^ny;q)_{n-k-1}$ for $k=0$ which may
happen only for $g\in G_n^{(1)}$. 
Suppose $g\notin G_n^{(1)}$. 
Then it is sufficient to use the fact that
any Laurent polynomial vanishing on a countable set vanishes identically.
It follows that, first, $f(1)\neq0$ for any non-zero L.-p.\ solution
and, second, any two non-zero L.-p.\ solutions are proportional, in 
particular to the standard solution $S^{(\ell;q)}_{\l}(y)$.

Instead of the sequence $y=q^{-j}$ one can take $y=q^j\ell^{-n}$ and
use the same argument. Note that the above recursive process is the exact
analog of the Taylor series
expansion around $y=1$ in the nonrelativistic case.

On the other hand, a similar argument works 
with expansion around $0$ or $\infty$.
Substituting in (\ref{eq:sep-eq-N}) the expansion
$f(y)=\sum_{k=k_-}^{k_+}f_ky^k$ one obtains $(n+2)$-terms recurrence
relation $\sum_{j=0}^{n+1}A_{kj}f_{k-j}=0$ for $f_k$. The ``boundary''
coefficients $A_{k0}$ and $A_{k,n+1}$ have simple form
\begin{subeqnarray}
A_{k0}&=&(-1)^n\ell^{\frac{n(n-1)}{2}}\prod_{j=1}^n(q^k-q^{\l_j}\ell^{1-j}),
\\
 A_{k,n+1}&=&-\left(\frac{\ell}{q}\right)^{\frac{n(n+1)}{2}}
\prod_{j=1}^n(q^k-q^{\l_j+n+1}\ell^{n-j}).
\label{eq:def-Ak}
\end{subeqnarray}

Suppose $g\notin G_n^{(0)}$. 
Then, since $\ell=q^{-g}$, the coefficient
$A_{k0}$ vanishes only for $k=\l_1$, and  $A_{k,n+1}=0$ only for
$k=\l_n+n+1$. Hence inevitably $k_-=\l_1$, $k_+=\l_n$, and the coefficients
$f_k$ are determined recursively in a unique way starting from $f_{\l_1}$
or $f_{\l_n}$ which proves the uniqueness of L.-p.\ solution.
\endproof

The question whether the uniqueness of the L.-p.\ solution really breaks for
$g\in G_n$, remains still open.

It would be interesting to strengthen the above result.

\begin{guess}
 The equation (\ref{eq:sep-eq-N}) with free parameters $h_j$
has a polynomial solution only for $h_j$ given by (\ref{eq:def-hk-gen})
and $\l=\{\l_1\leq \l_2\leq\ldots\leq \l_n\}\in\Z^n$.
\end{guess}

\newsection{Discussion}
The results of the present paper generalize to the case of the
Ruijsenaars model and Macdonald  polynomials those of \cite{Kuz-Skl}
obtained for the Calogero-Sutherland model and Jack polynomials.
In the nonrelativistic limit $\hbar\rightarrow0$, $g={\rm const}$,
the Hamiltonians $H_k$, operator $M$, separated polynomials $S$
and equations for them go over into the corresponding objects
described in \cite{Kuz-Skl}.

The crucial element of our approach is the operator identity
(\ref{eq:char-eq-M-op}) which allows to prove the factorization
(\ref{eq:act-M-Macd-3}) of $MP_\l^{(\ell;q)}$ and to establish thus
the separation of variables. The identity
(\ref{eq:char-eq-M-op}) is apparently a quantum analog of the
characteristic equation for the classical Lax operator.
Moreover, the kernel $\M$ can be considered as a collection of
eigenfunctions to the quantized separation variables $y_j$ describing thus
the change of basis from `$t$-representation' to `$y$-representation'.
Though these analogies with the classical inverse scattering
method proved to be
useful as an heuristic tool for finding SoV for quantum integrable
systems \cite{Skl-rev}, their algebraic/geometric origin
is still to be cleared up.

An interesting problem is to search for alternative forms of $M$.
We have presented here
two descriptions of $M$: analytical (\ref{eq:def-M}) in terms of
Askey-Wilson integral, and algebraic (\ref{eq:act-M-basis-t})  in terms of
the basis
$p_{jk\nu}$. Our study of $M$ is based mainly on the analytical definition.
It would be interesting also to develop the theory of $M$ based
entirely on the algebraic definition, in particular, to give a
purely algebraic proof of the identity  (\ref{eq:char-eq-M-op}).

When our work was close to be finished  we became aware of the preprint
\cite{Mangaz} of Mangazeev addressing the same problem of SoV for $A_2$
Macdonald polynomials. His proposal for the operator $M$ is different from
ours, using a $q$-integral rather than a contour integral as we do.
Some of his arguments are quite formal, for instance, expressions
with the ${}_6\psi{}_6$-series he is using as a final result are
divergent.
It seems that our choice of $M$, compared to that of \cite{Mangaz},
allows to overcome the problems of convergence of the $q$-integral
and to obtain explicit expressions for $M^{-1}$ and action of $M$ on
polynomials. Still, the problem of representing $M$ as a
$q$-integral seems to deserve a further consideration.

Although we can predict the form
of the separation polynomial $S_{\l}^{(\ell;q)}(y)$ for the $n$-particle case
and study it in detail (section 5),
the corresponding $n$-particle generalization
of the kernel ${\cal M}$ is not yet clear, so it is an open problem to
separate variables for the $A_{n-1}$ Macdonald polynomials for $n>3$.

In fact, there are infinitely many ``separating'' operators $M^{(n)}$,
since for any choice of $c_{\l}$ the operator defined as
\beq
 M^{(n)}: P_{\l}^{(\ell;q)}(t_1,\ldots,t_n)\rightarrow
         c_{\l}x^{\abs{\l}}
\prod_{j=1}^{n-1}S_{\l}^{(\ell;q)}(y_j)
\label{eq:act-M-Macd-N}
\eeq
will serve the purpose. The genuine problem, however, is to choose the
coefficients $c_{\l}$ in such a way that  the corresponding kernel
$\M^{(n)}$ were given by an explicit expression generalizing
(\ref{eq:def-ker-M}).

\setcounter{section}{0}\renewcommand{\thesection}{\Alph{section}}
\newappendix{Formulas from $q$-analysis}
 For reader's convenience, we have collected here the most important
definitions and formulas
from $q$-analysis used in the main text. For references see 
\cite{GR,Koorn-tut,Exton,Erdelyi}. 
Especially useful for practical calculations is the collection of formulas 
in Appendices I and II from \cite{GR}. Throughout the text it is assumed that
$0<q<1$.

The $q$-shifted factorial and its generalizations are defined as
\beq
 (a;q)_0:=1, \qquad 
 (a;q)_k:=(1-a)(1-aq)\ldots(1-aq^{k-1}), \qquad k=1,2,\ldots\ ,
\label{Poch-k}
\eeq
\beq
 (a;q)_\infty:=\prod_{k=0}^\infty(1-aq^k), \qquad
  (x;q)_\a=\frac{(x;q)_\infty}{(q^\a x;q)_\infty}, \qquad \a\in\C,
\label{eq:Poch-inf}
\eeq
\beq
 (a_1,a_2,\cdots,a_n;q)_k:=(a_1;q)_k(a_2;q)_k\ldots(a_n;q)_k,
\qquad k=0,1,2,\ldots \ {\rm or}\ \infty.
\label{eq:Poch-prod}
\eeq

Note the useful relations
\beq
 (qx;q)_\a=\frac{1-q^\a x}{1-x}(x;q)_\a, \qquad 
 (q^{-1}x;q)_\a=\frac{1-q^{-1}x}{1-q^{\a-1}x}(x;q)_\a.
\label{eq:shift-qx}
\eeq

We make use also of the $q$-binomial coefficient
\beq
 \left[{n \atop k}\right]_q:=\frac{(q;q)_n}{(q;q)_k(q;q)_{n-k}}, \qquad
k=0,1,\ldots,n,
\label{eq:def-q-binom}
\eeq
$q$-Gamma and $q$-Beta functions
\beq
 \Gamma_q(z)=\frac{(q;q)_\infty}{(q^z;q)_\infty(1-q)^{z-1}}, \qquad
 \Gamma_q(z+1)=\frac{1-q^z}{1-q}\Gamma_q(z),
\label{eq:def-q-Gamma}
\eeq
\beq
 B_q(a,b)=\frac{\Gamma_q(a)\Gamma_q(b)}{\Gamma_q(a+b)}=
    (1-q)\frac{(q,q^{a+b};q)_\infty}{(q^a,q^b;q)_\infty},
\label{eq:def-q-Beta}
\eeq
$q$-integral
\beq
 \int_0^1d_qt f(t):=(1-q)\sum_{k=0}^\infty f(q^k)q^k,
\label{eq:def-q-int}
\eeq
and basic hypergeometric series ($n\in\Z_{\geq0}$)
\beq
 {}_n\oldphi_{n-1}\left[
    \begin{array}{c} a_1,\ldots,a_n \\ b_1,\ldots,b_{n-1}\end{array};
     q,y\right]:=
\sum_{k=0}^\infty
   \frac{(a_1,\ldots,a_n;q)_k}{(q,b_1,\ldots,b_{n-1};q)_k}y^k,
\qquad \abs{y}<1.
\label{eq:def-bas-hgs}
\eeq

Denoting the expression (\ref{eq:def-bas-hgs}) by $f(y)$ we observe that it
satisfies the $n$-th order $q$-difference equation, see
\cite{Jackson} and \cite{Exton} (section 2.12.3):
\beq
\left\{
  y\prod_{k=1}^n(1-a_k Y)-\prod_{k=1}^n(1-q^{-1}b_k Y) 
\right\}f(y)=0
\label{eq:hg-diff-eq}
\eeq
where $(Yf)(y):=f(qy)$ and $b_n\equiv q$.

Summation formula for ${}_1\oldphi_0$ ($q$-binomial series):
\beq
 {}_1\oldphi_0\left[\begin{array}{c} a \\ - \end{array}; q,y\right]
\equiv \sum_{k=0}^\infty\frac{(a;q)_k}{(q;q)_k}y^k=
  \frac{(ay;q)_\infty}{(y;q)_\infty}, \qquad \abs{y}<1.
\label{eq:1phi0}
\eeq

The $\oldphi_D$-type $q$-Lauricella function
\cite{Erdelyi,q-Lauri} of $n-1$ variables $x_j$
is a multi-variable generalization of the basic hypergeometric series:
\beq
\oldphi_D\left[\begin{array}{l}
a';b'_1,\ldots,b'_{n-1}\cr c\end{array};q;x_1,\ldots,x_{n-1}\right]:=
\sum_{k_1,\ldots,k_{n-1}=0}^\infty
\frac{(a';q)_{k_1+\ldots+k_{n-1}}}{(c;q)_{k_1+\ldots+k_{n-1}}}
\prod_{j=1}^{n-1}\frac{(b'_j;q)_{k_j}x_j^{k_j}}{(q;q)_{k_j}}\,.
\label{eq:def-q-Lauri}
\eeq 

Andrews \cite{A72} has found that $\oldphi_D$ can be expressed in terms of
the basic hypergeometric function ${}_n\oldphi{}_{n-1}$ of one variable:
\be
\oldphi_D\left[\begin{array}{l}
a';b'_1,\ldots,b'_{n-1}\cr c\end{array};q;x_1,\ldots,x_{n-1}\right]=
\frac{(a',b'_1x_1,\ldots,b'_{n-1}x_{n-1};q)_\infty}{
(c,x_1,\ldots,x_{n-1};q)_\infty}
\ee
\beq \times
{}_n\oldphi{}_{n-1}\left[\matrix{c/a',&x_1,&\ldots ,&x_{n-1}\cr
{}&b'_1x_1,&\ldots,&b'_{n-1}x_{n-1}};q,a'\right]\,.
\label{eq:Andrews-Lauri}
\eeq

Our main technical tool, on which the proof of the main theorem 
\ref{th:main-th} depends, is the famous
Askey-Wilson integral (\cite{GR}, section 6.1; 
\cite{Koorn-tut}, section 2.6). Let
\beq
 w(a,b,c,d;t):=
\frac{(t^2,t^{-2};q)_\infty}%
{(at,at^{-1},bt,bt^{-1},ct,ct^{-1},dt,dt^{-1};q)_\infty}.
\label{eq:def-w-abcd}
\eeq

Then
\beq
 \frac{1}{2\pi i}\int\limits_{\G_{abcd}}\frac{dt}{t}\,w(a,b,c,d;t)=
\frac{2(abcd;q)_\infty}{(q,ab,ac,ad,bc,bd,cd;q)_\infty}.
\label{eq:AW-int}
\eeq

The cycle $\G_{abcd}$ depends on parameters $a$, $b$, $c$, $d$ and is 
defined as follows. Let $C_{z,r}$ be the counter-clockwise oriented circle 
with the center $z$ and radius $r$. 

If $\abs{a},\abs{b},\abs{c},\abs{d}<1$ then $\G_{abcd}=C_{0,1}$.
The identity (\ref{eq:AW-int})
can be continued analytically for the values of parameters $a$, $b$, $c$,
$d$ outside the unit circle provided the cycle $\G_{abcd}$ is deformed 
appropriately. In general case 
\beq
 \G_{abcd}=C_{0,1}+\sum_{x=a,b,c,d}
\sum_{\scriptstyle k\geq0\atop \scriptstyle \abs{x}q^k>1}
(C_{xq^k,\eps}-C_{x^{-1}q^{-k},\eps}),
\label{eq:def-contour}
\eeq
$\eps$ being small enough for $C_{x^{\pm1}q^{\pm k},\eps}$ to encircle 
only one pole of the denominator.

The following formulas are useful when studying the classical and 
nonrelativistic limits of the quantum Ruijsenaars model. Both correspond
to $\hbar\rightarrow0$, $q=e^{-\hbar}\rightarrow1$ and differ only in the
behaviour of $\ell$ (\ref{eq:def-ell-g}). As $q\uparrow1$,
\beq
 (x;q)_\a\rightarrow (1-x)^\a,
\eeq
\beq
 \frac{(q^\a;q)_k}{(1-q)^k}\rightarrow (\a)_k:=
\a(\a+1)\ldots(\a+k-1),
\eeq
\beq
{}_n\oldphi_{n-1}\left[
    \begin{array}{c} q^{\a_1},\ldots,q^{\a_n} \\ 
                     q^{\b_1},\ldots,q^{\b_{n-1}}\end{array};q,y\right]
\rightarrow
{}_nF_{n-1}\left[
    \begin{array}{c} \a_1,\ldots,\a_n \\ \b_1,\ldots,\b_{n-1}\end{array};
     y\right]
\eeq
where ${}_nF_{n-1}$ is the standard (generalized) hypergeometric series,
and finally (see \cite{dilog}, \mathhexbox278 2.5, Corollary 10):
\beq
 \ln(x;q)_\infty=-\hbar^{-1}\dilog(x)+\half\ln(1-x)+O(\hbar), \qquad
x\in(0,1).
\label{eq:asympt-dilog}
\eeq

\newappendix{Operator $M_{\a\b}$}

In this section the results are collected concerning the
two-parametric generalization $M_{\a\b}$ of one-parametric operator
family $M\equiv M_{g,2g}$ studied in the main text. It is an open question
whether $M_{\a\b}$ provides a SoV for some integrable model.

Let us substitute into the Askey-Wilson integral weight $w(a,b,c,d;t)$
(\ref{eq:def-w-abcd}) the values
\beq
 a=sq^{\frac\a2}, \qquad b=s^{-1}q^{\frac\a2}, \qquad
 c=rq^{\frac\b2}, \qquad d=r^{-1}q^{\frac\b2}, 
\label{eq:def-abcd}
\eeq
and introduce the notation (quantum analog of (\ref{eq:def-LL}))
\beq
 \Lq\nu xy:=(\nu xy,\nu xy^{-1},\nu x^{-1}y,\nu x^{-1}y^{-1};q)_\infty.
\label{eq:def-Lq}
\eeq

The kernel
\beq
 \K_{\a\b}(r,s\mid t):=
w(sq^{\frac\a2},s^{-1}q^{\frac\a2},rq^{\frac\b2},r^{-1}q^{\frac\b2};t)
=\frac{(t^2,t^{-2};q)_\infty}%
{\Lq{q^{\frac\a2}}st\Lq{q^{\frac\b2}}rt   }
\label{eq:def-kernel-Kab}
\eeq
defines the integral operator 
\beq
 K_{\a\b}:f(t)\rightarrow \frac{1}{2\pi i}\int\limits_{\G_{\a\b}^{rs}}
\frac{dt}{t}\,
\K_{\a\b}(r,s\mid t)f(t),
\label{eq:def-Kab}
\eeq
the contour $\G_{\a\b}^{rs}$ being obtained from $\G_{abcd}$ 
(\ref{eq:def-contour}) by substitutions (\ref{eq:def-abcd}).

Using the Askey-Wilson integral (\ref{eq:AW-int}) we obtain then the formula
\beq
 K_{\a\b}:1\rightarrow
\frac{2B_q(\a,\b)}{(1-q)(q;q)_\infty^2\,\Lq{q^{\frac{\a+\b}{2}}}rs  }.
\label{eq:act-K-1}
\eeq

Now we introduce the operator $M_{\a\b}$
\beq
 M_{\a\b}=(K_{\a\b}\cdot 1)^{-1}\circ K_{\a\b},
\label{eq:def-Mab}
\eeq
so that $M:1\rightarrow1$, and having the kernel
\beq
 \M_{\a\b}(r,s\mid t)=
\frac{(1-q)(q;q)_\infty^2(t^2,t^{-2};q)_\infty
\,\Lq{q^{\frac{\a+\b}{2}}}rs }%
{2B_q(\a,\b)
\,\Lq{q^{\frac\a2}}st\Lq{q^{\frac\b2}}rt  }.
\label{eq:def-ker-Mab}
\eeq

The kernel $\M$ (\ref{eq:def-ker-M}) studied in Section 4 is
obtained from $\M_{\a\b}$ (\ref{eq:def-ker-Mab}) 
after the substitutions
\beq
\a=g,\quad \b=2g,\quad q^{-g}=\ell, \quad
r=t_+=y_+\ell^{\frac32},\quad  s=y_-,\quad t=t_-.
\label{eq:ab-g}
\eeq

It is natural to think of  $M_{\a\b}$ as acting on the space
${\rm Ref}(t)$  of reflexive (invariant w.r.t.\ $t\rightarrow t^{-1}$)
Laurent polynomials in $t$. Consider a Laurent polynomial 
$R^{\a\b}_{j_1j_2k_1k_2}\in{\rm Ref}(t)$, $j_{1,2},k_{1,2}\in\Z_{\geq0}$
\begin{eqnarray}
 R^{\a\b}_{j_1j_2k_1k_2}(t)&:=&
   (q^{\frac{\a}{2}}st,q^{\frac{\a}{2}}st^{-1};q)_{j_1}
(q^{\frac{\a}{2}}s^{-1}t,q^{\frac{\a}{2}}s^{-1}t^{-1};q)_{j_2}
\nonumber \\
 &\times&(q^{\frac{\b}{2}}rt,q^{\frac{\b}{2}}rt^{-1};q)_{k_1}
(q^{\frac{\b}{2}}r^{-1}t,q^{\frac{\b}{2}}r^{-1}t^{-1};q)_{k_2}.
\label{eq:def-poly-P}
\end{eqnarray}

Using the obvious identity
\beq
 \K_{\a\b}(r,s\mid t)R^{\a\b}_{j_1j_2k_1k_2}(t)=
\K_{\a+j_1+j_2,\b+k_1+k_2}
(rq^{\frac{k_1-k_2}{2}},sq^{\frac{j_1-j_2}{2}}\mid t),
\label{eq:KP}
\eeq
together with (\ref{eq:act-K-1}),
we obtain the formula for action of $M_{\a\b}$ on Laurent polynomials 
\begin{eqnarray}
M_{\a\b}: R^{\a\b}_{j_1j_2k_1k_2}&\rightarrow&
\frac{(q^\a;q)_{j_1+j_2}(q^\b;q)_{k_1+k_2}}%
{(q^{\a+\b};q)_{j_1+j_2+k_1+k_2}} 
\nonumber \\
&&\times
(q^{\frac{\a+\b}{2}}rs;q)_{j_1+k_1}
(q^{\frac{\a+\b}{2}}rs^{-1};q)_{j_2+k_1} \nonumber \\
&&\times
(q^{\frac{\a+\b}{2}}r^{-1}s;q)_{j_1+k_2}
(q^{\frac{\a+\b}{2}}r^{-1}s^{-1};q)_{j_2+k_2}.
\label{eq:act-Mab-P}
\end{eqnarray}

The set of polynomials $R^{\a\b}_{j_1j_2k_1k_2}$ is rich enough 
to choose from it a basis in ${\rm Ref}(t)$, for instance
\beq
p_\nu^\b(t):=R^{\a\b}_{00\nu0}(t)
\equiv(q^{\frac\b2}rt,q^{\frac\b2}rt^{-1})_\nu,
\qquad \nu=0,1,2,\ldots\,.
\label{eq:def-pm}
\eeq

More correctly,
since $p_\nu^\b(t)=(-1)^\nu q^{\nu(\nu-1+\b)/2}r^\nu(t^\nu+t^{-\nu})+%
${\it lower order terms}, $p_\nu^\b(t)$ is a basis in the space of
reflexive Laurent polynomials in variable $t$ with coefficients from
${\rm Ref}(r)$. The specialization of the formula (\ref{eq:act-Mab-P})
\beq
 M_{\a\b}:p_\nu^\b(t)\rightarrow\frac{(q^\b;q)_\nu}{(q^{\a+\b};q)_\nu}
p_\nu^{\a+\b}(s)
\label{eq:act-Mab-pm}
\eeq
provides thus a tool for calculating explicitely the action of $M$
on any polynomial $\in{\rm Ref(t)}$.

Analysing the formula (\ref{eq:act-Mab-pm}) one obtains the following
statement.

\begin{prop}
Let $f\in{\rm Ref}(t)$ and suppose $f$ has the parity $\sigma$ that is
$f(-t)=(-1)^\sigma f(t)$. Let $M_{\a\b}:f\rightarrow F$. Then:
$F\in{\rm Ref}(r)\otimes{\rm Ref}(s)$,
$F(r,s)=\left.F(s,r)\right|_{\a\leftrightarrow\b}$,
$F(-r,-s)=(-1)^\sigma F(r,s)$.
\label{pr:act-Mab-poly}
\end{prop}

The integral operator $M_{\a\b}$ simplifies drastically when one of the
parameters $\a$, $\b$ takes negative integer values.
\begin{prop}
Let $\a\in\Z_{\leq0}$. Then $M_{\a\b}$ turns into the 
$q$-difference operator of order $-\a$:
\beq
 M_{\a\b}: f(t)\rightarrow \sum_{k=0}^{-\a} \xi_k(r,s)f(q^{k+\frac{\a}{2}}s)
\label{eq:Mab-fin-dif}
\eeq
where
\begin{eqnarray}
\lefteqn{
 \xi_k(r,s)=(-1)^k q^{-\frac{k(k-1)}{2}} 
\left[\begin{array}{c}-\a \\ k\end{array}\right]_q s^{-2k}
(1-q^{-\a-2k}s^{-2}) }\nonumber \\
&&\times
\frac{(q^{\frac{\a+\b}{2}}rs,q^{\frac{\a+\b}{2}}r^{-1}s;q)_k
(q^{\frac{\a+\b}{2}}rs^{-1},q^{\frac{\a+\b}{2}}r^{-1}s^{-1};q)_{-\a-k}}%
{(q^{\a+\b};q)_{-\a}(q^{-k}s^{-2};q)_{1-\a} }.
\label{eq:expr-xik}
\end{eqnarray}
\label{pr:Mab-int}
\end{prop}

{\bf Proof.} Instead of analyzing the degeneration of the integral operator
defined by the kernel (\ref{eq:def-ker-Mab}) it is easier to
study the action of $M_{\a\b}$ on the basic
polynomials $p_\nu^\b(t)$ (\ref{eq:def-pm}).

Substituting $f(t):=p_\nu^\b(t)$ into (\ref{eq:Mab-fin-dif})
and using (\ref{eq:act-Mab-pm}) we obtain,
after simplification, the equality
\beq
 \sum_{k=0}^{-\a} \xi_k(r,s)
\frac{(q^{\frac{\a+\b}{2}+\nu}rs;q)_k
(q^{\frac{\a+\b}{2}+\nu}rs^{-1};q)_{-\a-k}}%
{(q^{\frac{\a+\b}{2}}rs;q)_k(q^{\frac{\a+\b}{2}}rs^{-1};q)_{-\a-k}}=
\frac{(q^{\a+\b+\nu};q)_{-\a}}{(q^{\a+\b};q)_{-\a}}
\label{eq:eq-xik}
\eeq
which, after substituting (\ref{eq:expr-xik}) and making a series of
elementary transformations (see Appendix I to  \cite{GR}),
can be put into the form

\be
\sum_{k=0}^{-\a} q^{(1-\a-\b-\nu)k}
\frac{(q^{\a},q^{\frac{\a+\b}{2}+\nu}rs,q^{\frac{\a+\b}{2}}r^{-1}s,q^{\a}s^2,
q^{1+\frac{\a}{2}}s,-q^{1+\frac{\a}{2}}s;q)_k}%
{(q,q^{\frac{\a-\b}{2}+1}rs,q^{\frac{\a-\b}{2}+1-\nu}r^{-1}s,qs^2,
q^{\frac{\a}{2}}s,-q^{\frac{\a}{2}}s;q)_k}
\ee\beq
=\frac{(q^{1-\b-\nu},q^{1+\a}s^2;q)_{-\a}}%
{(q^{\frac{\a-\b}{2}+1}rs,q^{\frac{\a-\b}{2}+1-\nu}r^{-1}s;q)_{-\a}},
\eeq
identical to the summation formula (II.21) from \cite{GR}:
\beq
 {}_6\oldphi_5\left[\begin{array}{c}
 a,qa^{\half},-qa^{\half},b,c,q^{\a} \\
a^{\half},-a^{\half},\dfrac{aq}{b},\dfrac{aq}{c},aq^{1-\a}
\end{array};q,\frac{aq^{1-\a}}{bc}\right]=
\frac{\left(aq,\dfrac{aq}{bc};q\right)_{-\a}}%
{\left(\dfrac{aq}{b},\dfrac{aq}{c};q\right)_{-\a}}
\eeq
for the following identification of the parameters
\beq
 a=q^{\a}s^2, \quad b=q^{\frac{\a+\b}{2}}r^{-1}s, \quad
c=q^{\frac{\a+\b}{2}+\nu}rs.
\eeq

By the symmetry $\a\leftrightarrow\b$, $r\leftrightarrow s$ the similar
statement can be proved for $\b\in\Z_{\leq0}$.
\endproof

To determine the inversion of $M_{\a\b}$
let us think of  $r$ as a parameter and of
$M_{\a\b}$ as an operator 
$M_{\a\b}^r:{\rm Ref}(t)\rightarrow{\rm Ref}(s):f(t)\rightarrow F(s)$.
Then,
applying $M_{\a\b}^r$ to the basis $p_\nu^\b(t)$ (\ref{eq:def-pm}) and 
inverting the formula (\ref{eq:act-Mab-pm}), we obtain the following
statement.

\begin{prop}
The inversion formula for $M_{\a\b}^r$:
\beq
 \left(M_{\a\b}^r\right)^{-1}=M_{-\a,\a+\b}^r.
\label{eq:inv-Mab}
\eeq
The corresponding kernel is
\beq
 \tilde\M_{\a\b}^r(t\mid s)=
\frac{(1-q)(q;q)_\infty^2(s^2,s^{-2};q)_\infty
\,\Lq{q^{\frac{\b}{2}}}rt  }%
{2B_q(-\a,\a+\b)
\,\Lq{q^{-\frac{\a}{2}}}st\Lq{q^{\frac{\a+\b}{2}}}rs   }.
\label{eq:def-ker-Mab-inv}
\eeq
\label{pr:inv-Mab}
\end{prop}
\pagebreak

\end{document}